\documentclass[english,10pt]{article}
\usepackage{enumitem}
\usepackage{cite}
\usepackage{ulem}
\usepackage{amsthm}
\usepackage[latin9]{inputenc}
\usepackage{babel}
\usepackage[lmargin=1.0in, rmargin=1.0in, tmargin=0.75in, bmargin=0.75in]{geometry}
\usepackage{float}
\usepackage{subfig}
\usepackage{graphicx}
\usepackage{caption}
 \usepackage[colorlinks=true,citecolor=blue,urlcolor=blue]{hyperref}
\usepackage{multirow}
\usepackage[title]{appendix}
\usepackage{url}
\usepackage{bm}
\usepackage[section]{placeins}
\usepackage{placeins}
\usepackage{amsmath}
\usepackage{amssymb,amsfonts}
\usepackage{upgreek}
\usepackage{lscape}
\usepackage{rotating}
\usepackage{comment}
\usepackage{physics}
\usepackage{color}
\usepackage[numbers]{natbib}

\numberwithin{theorem}{subsection} 
\newtheorem{remark}{Remark}

\numberwithin{remark}{subsection}

\usepackage{verbatim}
\usepackage{listings}

\definecolor{mygreen}{rgb}{0,0.6,0}
\definecolor{mygray}{rgb}{0.5,0.5,0.5}
\definecolor{mymauve}{rgb}{0.58,0,0.82}

\usepackage{MyMacros}



\begin{document}
\title{An unstructured finite element model for incompressible
  two-phase flow based on a monolithic conservative level set method}
\author{
  Manuel Quezada de Luna
  \thanks{KAUST CEMSE, manuel.quezada@kaust.edu.sa. Research performed as a post-doc at USACE-ERDC-CHL.} \and
  J. Haydel Collins
  \thanks{USACE-MVN, jason.h.collins@usace.army.mil. Research performed as a participate in USACE-ERDC University.} \and
  Christopher E. Kees
  \thanks{USACE-ERDC-CHL, christopher.e.kees@usace.army.mil.} 
}

\date{}
\maketitle


\begin{abstract}
We present a robust numerical method for solving 
incompressible, immiscible two-phase flows.
The method extends the monolithic phase conservative
level set method with embedded redistancing by \citet{quezada2018monolithic}
and a semi-implicit high-order projection scheme for variable-density flows
by \citet{guermond2009splitting}.
The level set method can be initialized conveniently via a simple phase indicator field,
which is pre-processed to obtain an approximate signed distance function.
To do this, we propose a new PDE-based redistancing method.
We also improve the scheme in \cite{quezada2018monolithic} to provide more accuracy and robustness in full two-phase flow simulations. Specifically, we perform an extra step to ensure convergence to the signed distance level set function and simplify other aspects of the original scheme.
Lastly, we introduce consistent artificial viscosity to stabilize the momentum equations in the context of the projection scheme.
This stabilization is algebraic, has no tunable parameters and
is suitable for unstructured meshes and arbitrary refinement levels.
The overall methodology includes few numerical tuning parameters;
however, for the wide range of problems that we solve, we identify
only one parameter that strongly affects performance of the
computational model and provide a value that provides accurate results across all the benchmarks presented.
The result is a robust, accurate, and efficient two-phase flow model, which is mass- and volume-conserving on unstructured meshes and has low user input requirements for real applications.
\end{abstract}

\section{Introduction}

Understanding the behavior of immiscible fluids with distinct material characteristics
(e.g. water and air) is important for many engineering and industrial applications.
Water-oil-gas interaction within subterranean reservoirs, combustion engines, and open channel
hydraulics are all examples of flow phenomena involving multiple fluids, i.e. multiphase
flow. In order to accurately describe the evolving interface between individual fluid
subdomains, high fidelity numerical models are required. The methods involved should aim
for computational efficiency along with preservation of qualitative properties of the
material continuum. In the case of multiphase simulations concerning divergence-free
velocity fields, maintaining the conservation of mass and volume is imperative.

In general, multiphase flow is modeled by solving the Navier-Stokes equations with
spatially varying material parameters.
We consider an Eulerian description for the fluid motion along with continuous Galerkin finite elements
for discretization of all equations in space. For Eulerian methods, the fluid motion is
calculated from a stationary frame of reference, a fixed computational grid. These
techniques can resolve large fluid deformations, however, a representation of
the phase interface is required.

Two routinely implemented methods for representing and evolving
interfaces between discrete phases are the
Volume of Fluid (VOF) method by \cite{hirt1981volume} 
and level set techniques by \cite{osher1988fronts, sussman1994level}.
The VOF method assigns phase identities, e.g. fluids A and B,
via a characteristic function that equals one in fluid A
and zero in fluid B. Grid cell averaging of this phase indicator function results in
the designation of a cell volume fraction in the range [0, 1]. An interface is then reconstructed in cells with intermediate volume fraction and the VOF re-initialized. For many multiphase flow applications,
the characteristic function is transported using 
velocities obtained from the (incompressible) Navier-Stokes equations.
This solution produces density and viscosity fields requiring subsequent
interfacial reconstruction. Once the phase boundary location is accurately
reconstructed, material properties can then be designated respectively for each
fluid subdomain.

In the level set method, the interface between the phases is represented
implicitly by a level set of a scalar function defined on the entire domain; e.g., as the zero level set of the signed distance function (SDF)
to the interface. The SDF returns a zero value at the interface and either a positive
or a negative distance for subdomains fluid A and fluid B respectively.
The level set function is advected by the fluid velocity field produced by the
(incompressible) Navier-Stokes equations.

Both the VOF and the level set methods bear limitations when
implemented alone.  Despite maintaining volume (or phase)
conservation, the VOF requires reconstruction of interface locations
from less precise cell averages.  The interface is effectively smeared
over some distance specific to the numerical scheme and then
reconstructed to be sharp base on geometric approximations. Level set
methods do not require interface reconstruction; however, they lack a
discrete conservation property; phases enclosed by the interface can
suffer from a significant loss of volume due to the accumulation of
discrete conservation errors over time.

There are hybrid methods combining ideas of the VOF, level set, and
particle methods. See for instance \cite{jcph.2000.6537, Ianniello,
  Enright_Hybrid_PLS_Method}.  In this work we consider a new hybrid
level set and volume-of-fluid that builds on the weak variational form
of phase conservation proposed in \cite{kees2011conservative}, the optimal control formulation of \cite{basting2014optimal}, and the
recent reformulation and extension of these approaches to a monolithic
scheme in \cite{quezada2018monolithic}.  These methods simultaneously maintain the convenient and precise signed distance represenation of the phase geometry and a phase conservation property.

Despite the numerical convenience of the signed distance field, it can be difficult for model users to provide initial conditions, since it generally requires solution of the nonlinear Eikonal equation from some other description of the interface.  In this work we provide a method that only requires input of the initial phase indicator field, which we pre-process to obtain the SDF. To do this, we solve the
non-linear Eikonal equation with a stabilization term that blends
ideas from \cite{lohmann2017flux}, the elliptic redistancing by
\cite{basting2014optimal} and \cite{quezada2018monolithic}.
In addition, we modify the method in \cite{quezada2018monolithic} by
simplifying the numerical discretization and solving (at every time
step) a redistancing pre-stage based on an extensions of the
stabilized Hamilton-Jacobi formulation of redistancing from
\cite{kees2011conservative}, which ensures convergence of the elliptic
redistancing approach to the correct SDF.

The velocity field is obtained via the second order projection 
scheme for the incompressible Navier-Stokes equations with 
variable material parameters by \cite{guermond2009splitting}.
We use Taylor-Hood finite elements. Additionally, we 
incorporate artificial viscosity based on \cite{FCT-book,badia2017monotonicity,barrenechea2016edge,guermond2017secondOrderdij}
and surface tension as proposed in \cite{hysing2006new}.
The artificial viscosity has no tunable parameters and is algebraic,
which makes it robust and suitable for unstructured meshes and arbitrary refinement levels.

The rest of this work is organized as follows.
In \S\ref{sec:finite_element_discretization}
we describe the finite element spatial discretization.
Afterwards, in \S\ref{sec:splitting_algorithm} we describe the overall algorithm.
Later, in \S\ref{sec:level_set_method} we review the conservative
level set method by \cite{quezada2018monolithic}. We consider only the
continuous model and propose a simpler discretization in time.
Furthermore, we propose a redistancing pre-stage that can be used
to obtain the SDF from a discontinuous indicator field.
In \S\ref{sec:navier-stokes} we describe the full Navier-Stokes discretization.
In this section we propose an algebraic and robust stabilization for the momentum equations.
Section \ref{sec:numerical_examples} is devoted to the numerical examples.
Finally, we close with some discussion in \S\ref{sec:conclusions}.
%

\section{Finite element spatial discretization}\label{sec:finite_element_discretization}
We use continuous Galerkin finite elements to discretize all equations in space.
Let $d=\{1,2,3\}$ be the number of spatial
dimensions and $\Omega\subset\mathbb{R}^d$ be a bounded domain with
boundary $\partial\Omega\subset\mathbb{R}^{d-1}$ which we decompose
into $\partial\Omega^-=\{\bfx\in\partial\Omega \st \bfu\cdot\bfn<0\}$
and $\partial\Omega^+=\{\bfx\in\partial\Omega \st \bfu\cdot\bfn\geq0\}$.
Given an SDF $\phi$ we define
$\Gamma(t)=\{\bfx\in\Omega \st \phi(\bfx,t)=0\}$
to be the material interface. 
Time-dependent variables are defined on the time interval $t\in [0,T]$, where $T>0$.
Given a computational mesh $\mathcal{T}_h$, we consider the finite element space
$X_h^p=\{ w\in\mathcal{C}^0(\Omega) \st w|_K\in\mathbb{P}_\text{p}, \forall K\in\mathcal{T}_h \}$
with $p=\{1,2\}$; i.e., we use only continuous piecewise linear and quadratic spaces. 
The spaces are spanned by basis functions
$\{w_1, \dots, w_{\text{dim}\left(X_h^p\right)}\}$ which possess
the partition of unity property; i.e., $\sum_jw_j(\bfx)=1$. The
degrees of freedom associated with these basis functions are
denoted by uppercase letters.
The finite element
solution $u_h(\bfx)\in X_h$ is given by
$u_h(\bfx)=\sum_{j\in\I(\Omega_i)}U_jw_j(\bfx)$, where $\Omega_i$
is the patch of elements containing node $i$. Here, and in the
rest of this paper, the notation $\I(z)$ is used for the index set
containing the numbers of all basis functions whose support on $z$
is of nonzero measure.

\section{Two-phase flow algorithm based on operator splitting}\label{sec:splitting_algorithm}
The overall two-phase flow algorithm is driven by a simple first-order operator splitting technique,
as in \cite{lin2005level,smolianski2005finite,kees2011conservative}. 
In particular, operator splitting is used to decouple the level set and the Navier-Stokes stages.
More sophisticated alternatives are possible as well, see for instance \cite{knio1999semi,nagrath2005computation,marchandise2006stabilized}. 
Following this approach one can discretize first the level set equation and subsequently the Navier-Stokes
equations or vice versa; we choose the former, which is analogous to the splitting presented for variable-density flows in \cite{guermond2009splitting}. 
In the rest of this section, for simplicity of exposition, we assume a fixed time step;
however, we apply this algorithm to a variable time stepping scheme, see \S\ref{sec:projection_scheme}.
Let $\phi$, $\bfu$ and $p$ denote the level set, velocity and pressure fields respectively. 
For any given time step assume we know the solution at time $t^n$ and $t^{n-1}$; in particular,
that we know $\phi^n$, $\bfu^n$, $\bfu^{n-1}$, $p^n$ and $p^{n-1}$.
We proceed as follows:
\begin{enumerate}
\item Obtain a second order extrapolation of the velocity field: $\bfu^*=2 \bfu^n - \bfu^{n-1}\approx\bfu^{n+1}$.
\item Using $\bfu^*$, solve the level set equation described in \S\ref{sec:level_set_method}, to obtain $\phi^{n+1}$.
\item Given the level set solution at $t^{n+1}$, define the fluid density ($\rho^{n+1}$) and dynamic viscosity ($\mu^{n+1}$) via
  \begin{align*}
    \rho^{n+1} &= \rho_A H_\epsilon(\phi^{n+1}) + \rho_W [1-H_\epsilon(\phi^{n+1})], \\
    \mu^{n+1} &= \mu_A H_\epsilon(\phi^{n+1}) + \mu_W [1-H_\epsilon(\phi^{n+1})],
  \end{align*}
  where $A$ and $W$ denote the air and water phases respectively, and $H_\epsilon$ is a regularized Heaviside function
  defined in \S\ref{sec:level_set_method}. 
\item Given $\rho^{n+1}$ and $\mu^{n+1}$, solve the variable-coefficient Navier-Stokes equations via a projection scheme,
  see \S\ref{sec:projection_scheme}, to obtain an updated velocity field $\bfu^{n+1}$.
\item Repeat until the final time is reached. 
\end{enumerate}

\section{Monolithic conservative level set method}\label{sec:level_set_method}

First note that a weak formulation of phase volume conservation for phase $\alpha$, with phase indicator function $\chi_{\alpha}$, and boundary described by the zero level set of an SDF $\phi$, can be written as
\begin{equation}\label{conservation}
  \int_{\Omega^*}  (H(\phi) - \chi_{\alpha}) w d\bfx = 0 \quad \forall w \in W(\Omega^*)  
\end{equation}
where $W(\Omega^*)$ is any suitable test space containing the constant function on subsets $\Omega^*$ of $\Omega$. For example, global conservation is ensured if $\Omega^*=\Omega$ or local conservation if $\Omega^*=K$. Equation \eqref{conservation} is ill-posed as an equation for $\phi$, but, if $W(\Omega^*)$ has the partition of unity property, then the conservation property is maintained even after addition of some regularization $G$ \cite{kees2011conservative}:
\begin{equation}\label{conservationReg}
  \int_{\Omega^*}  \left[ (H(\phi) - \chi_{\alpha}) w + G \nabla w \right] d\bfx = 0 \quad \forall w \in W(\Omega^*)  
\end{equation}

In this work we consider a time-dependent form of equation \eqref{conservationReg} from \cite{quezada2018monolithic} for solenoidal velocity field $\bfv$, which is given in strong form by
\begin{subequations}\label{level_set_model}
\begin{align}
  \partial_t \Seps(\phi)
  +\nabla\cdot\left[\bfv\Seps(\phi)-\lambda(\nabla\phi-\bfq)\right] &= 0, & \forall & \bfx\in\Omega
  \label{level_set_model_eqn1}\\
  \sqrt{|\nabla\phi|^2+\delta^2}\bfq &= \nabla\phi, & \forall & \bfx\in\Omega,
  \label{level_set_model_eqn2}\\
  \Seps(\phi) &= \Seps\left(\phi^\text{BC}\right),  & \forall & \bfx\in\partial\Omega^-, \\
  (\nabla\phi-\bfq)\cdot\bfn &=0, & \forall & \bfx\in\partial\Omega
\end{align}
\end{subequations}
where $\phi\in\mathbb{R}$ denotes the SDF level set function, $\Seps$
is a smoothed sign function, $\lambda>0$ is a user defined parameter
and $\delta>0$ is a small regularization parameter.
We use $\delta=1\times 10^{-10}$ in all simulations.
The smoothed sign function makes conservation symmetric with respect to the phases and
the regularization is based on an elliptic form of the Eikonal
equation\cite{basting2014optimal}. Equation \eqref{level_set_model}
combines level set evolution, signed distance property, and VOF
evolution into a single equation as described further below.
We follow \cite{quezada2018monolithic} and consider
\begin{align*}
    \Seps(\phi) &= 2\Heps(\phi)-1, \\
  \Heps(\phi)&=
  \begin{cases}
    0, & \text{ if } \phi\leq-\epsilon, \\
    \frac{1}{2}\left(1 + \frac{\phi}{\epsilon}
    + \frac{1}{\pi}\sin(\pi\frac{\phi}{\epsilon})
    \right),
    &\text{ if } -\epsilon < \phi < \epsilon, \\
    1, & \text{ if } \phi\geq\epsilon.
  \end{cases}
\end{align*}
Let $h(\bfx)$ denote a characterization of the local mesh size,
then $\epsilon=\alpha h(\bfx)$ defines the thickness of the regularization
of $\Heps$ and $\Seps$. We use $\alpha=\frac{3}{2}$ in all simulations. 

Equation \eqref{level_set_model_eqn1} is a conservation law for $\Seps(\phi)$. Under appropriate
boundary conditions and a conservative finite element method (e.g. the partition of unity property), this implies that
\begin{align*}
  \partial_t \int_\Omega \Seps\left(\phi(\bfx,t)\right) d\bfx = 0 \implies
  \partial_t \int_\Omega \Heps\left(\phi(\bfx,t)\right) d\bfx = 0,
\end{align*}
and since $\int_\Omega \Heps(\phi) d\bfx$ represents the (regularized) volume of one
of the phases, the method is volume conservative. 
The terms $\partial_t\Seps(\phi)+\nabla\cdot \bfv\Seps(\phi)=0$ correspond to a volume of fluid
like model that impose not only conservation but are also responsible for advecting
the interface to the correct position.
This is true since the velocity field is assumed to be solenoidal ($\nabla\cdot \bfv =0$)
and $\partial_t\Seps(\phi)+\bfv\cdot\nabla\Seps(\phi)=0 \implies \partial_t\phi+\bfv\cdot\nabla\phi=0$. 
Therefore, these terms are important for consistency with the volume of fluid
and level set methods. Henceforth, we refer to them as consistency terms of
\eqref{level_set_model_eqn1}.
Since the Jacobian of $\partial_t\Seps(\phi)+\nabla\cdot \bfv\Seps(\phi)=0$ is singular,
the terms $-\nabla\cdot\lambda(\nabla\phi-\bfq)$ are added to the model; i.e., they act as
a regularization.
Note that $\bfq\approx\frac{\nabla\phi}{|\nabla\phi|}$; hence,
the regularization terms also penalize deviations of the level set
from the distance function.
The result is a model for the SDF level set that is volume conservative and contains
a term that penalizes deviations from the distance function. 

\begin{remark}[Boundary conditions]
  Note that technically one has to impose boundary conditions for the level set $\phi$
  in the inflow boundary $\partial\Omega^-$. However, this information is only applied
  through the smoothed sign function $\Seps$. Therefore, all we need to know is
  $\Seps\left(\phi^\text{BC}\right)$; i.e., the phase (e.g., water or air)
  at that particular inflow boundary. 
  The second boundary condition $(\nabla\phi-\bfq)\cdot\bfn=0,~\forall\bfx\in\partial\Omega$,
  is imposed to guarantee conservation. Since the regularization term is integrated by parts during
  the finite element spatial discretization, applying this boundary condition is trivial. 
\end{remark}

\begin{remark}[About the parameter $\lambda$]\label{remark:about_lambda}
  The parameter $\lambda$ controls the amount of regularization and penalization introduced to
  the model.
  For any given problem, it is important to guarantee consistency of \eqref{level_set_model_eqn1}
  w.r.t. the volume of fluid equation. Therefore, the authors in \cite{quezada2018monolithic}
  propose to scale $\lambda$ by $h(\bfx)$. We follow their definition and use
  \begin{align*}
    \lambda=\tilde{\lambda}\left(\frac{v^{\max}}{C}\right)\frac{h(\bfx)}{|||\phi|-\overline{|\phi|}||_{L^\infty(\Omega)}},
  \end{align*}
  where $\tilde{\lambda}=\bigO(1)$ is a dimensionless user defined parameter,
  $v^{\max}=\max_\bfx |\bfv|$, $C<1$ is the Courant number and
  $\overline{|\phi|}=\frac{1}{|\Omega|}\int_\Omega|\phi| d\bfx$.

  We identify this parameter as the most important in the overall methodology. In all the problems
  in \S\ref{sec:numerical_examples} we obtain satisfactory results using $\tilde{\lambda}=10$. 
  However, it is essential to realize that the size of $|||\phi|-\overline{|\phi|}||_{L^\infty(\Omega)}$
  changes depending on the domain. In particular, if the domain allows the maximum value of $|\phi|$ to grow
  from one problem to another, then the effective value of $\lambda$ is reduced. If $\lambda$ is too small
  we observe perturbations in the free surface. And if $\lambda$ is too large more dissipation is added.
  In general, for any given problem, we want to choose $\tilde{\lambda}$ to be as small as possible without
  introducing large perturbations to the free surface.
  %
  In \cite{quezada2018monolithic}, the authors propose to use automated control theory to adjust this parameter at every time step.
  See for instance \cite{basting2014optimal}. 
\end{remark}

\subsection{Discontinuous initialization of the distance function level set}
\label{sec:disc_initialization}

\begin{figure}[!htb]
  \subfloat[$\hat{\phi}$]{\includegraphics[scale=0.19]{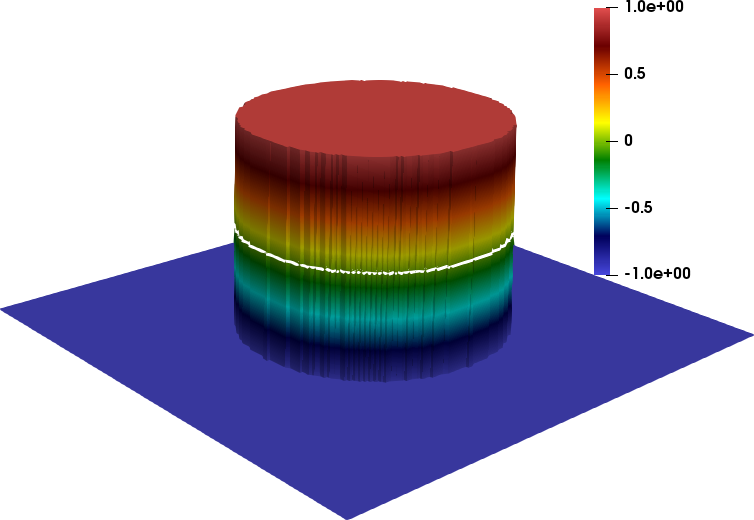}
    \label{fig:disc_ICs_phiHat}}\quad
\subfloat[Interface $\hat{\Gamma}$]{\includegraphics[scale=0.09]{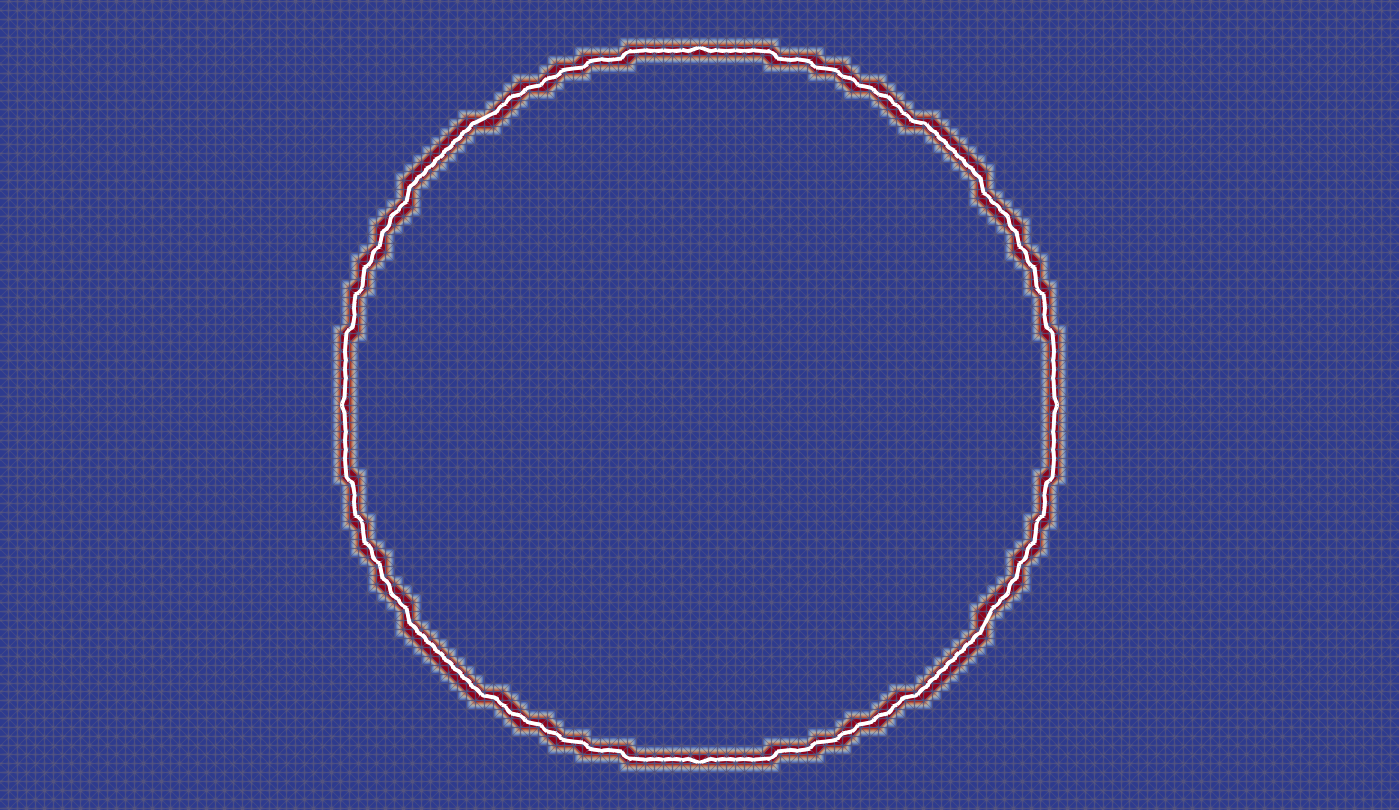}}\quad  
\subfloat[Step 1: $\hat{\phi}_h$]
         {\includegraphics[scale=0.19]{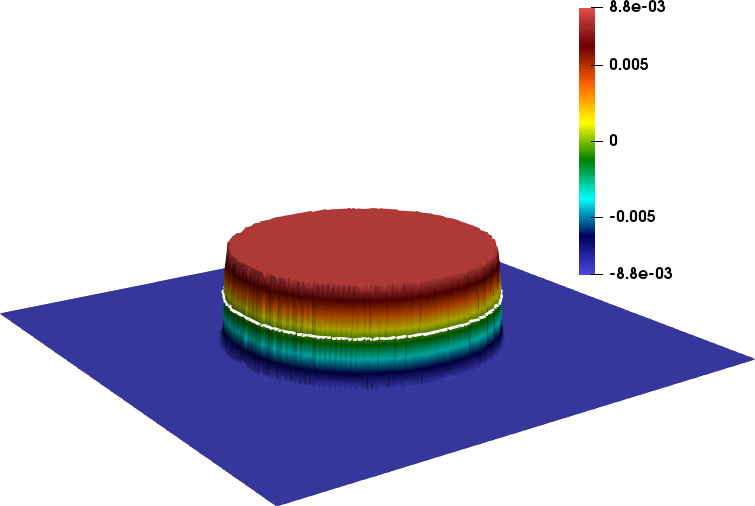}
         \label{fig:disc_ICs_phi_hHat}}\qquad

\subfloat[Step 2: $\phi^*_h$]{\includegraphics[scale=0.13]{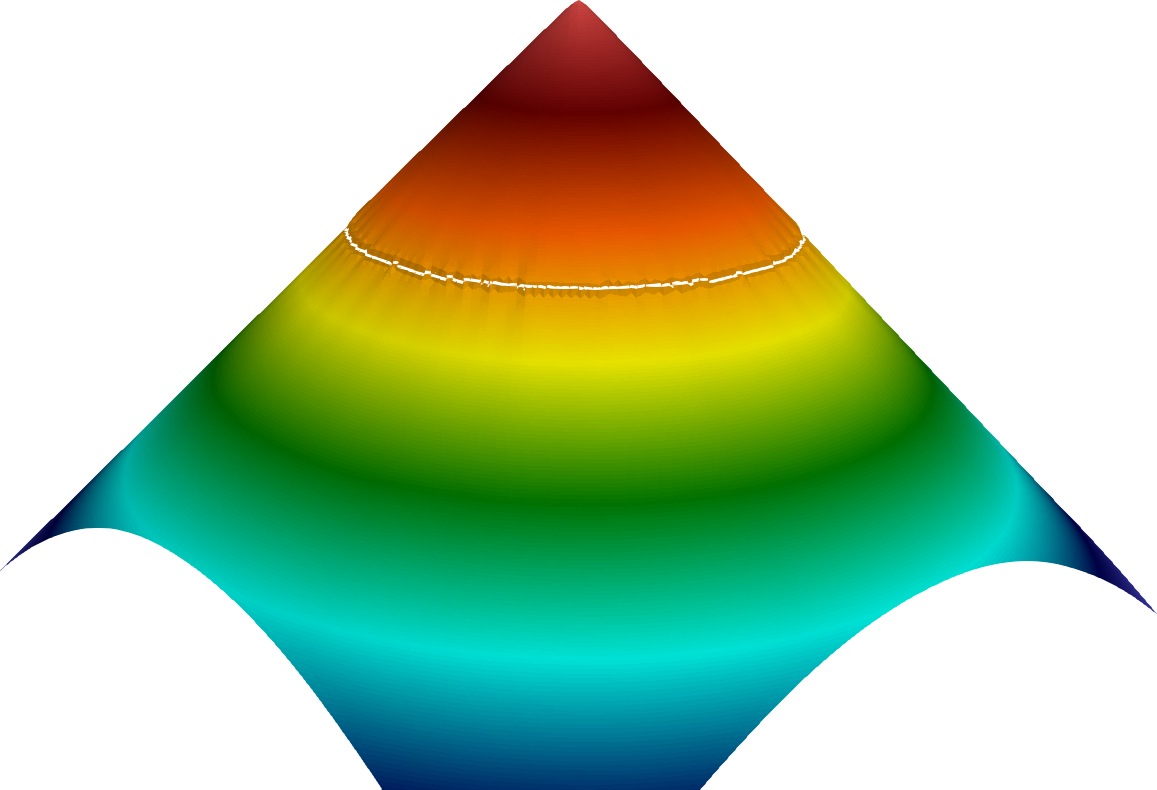}
  \label{fig:disc_ICs_phi_star}}\quad
\subfloat[Step 3: $\phi^{**}_h$]{\includegraphics[scale=0.13]{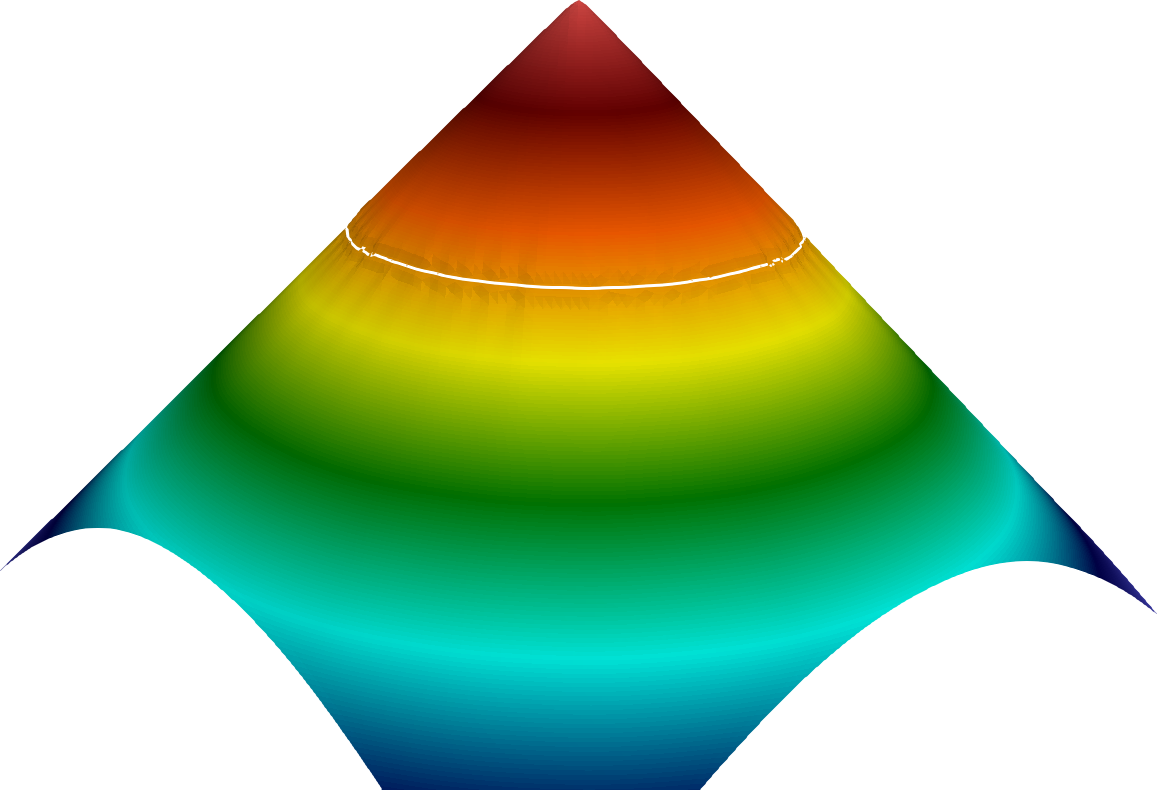}
  \label{fig:disc_ICs_phi_sstar}}\quad
\subfloat[Step 4: $\phi^0_h$]{\includegraphics[scale=0.13]{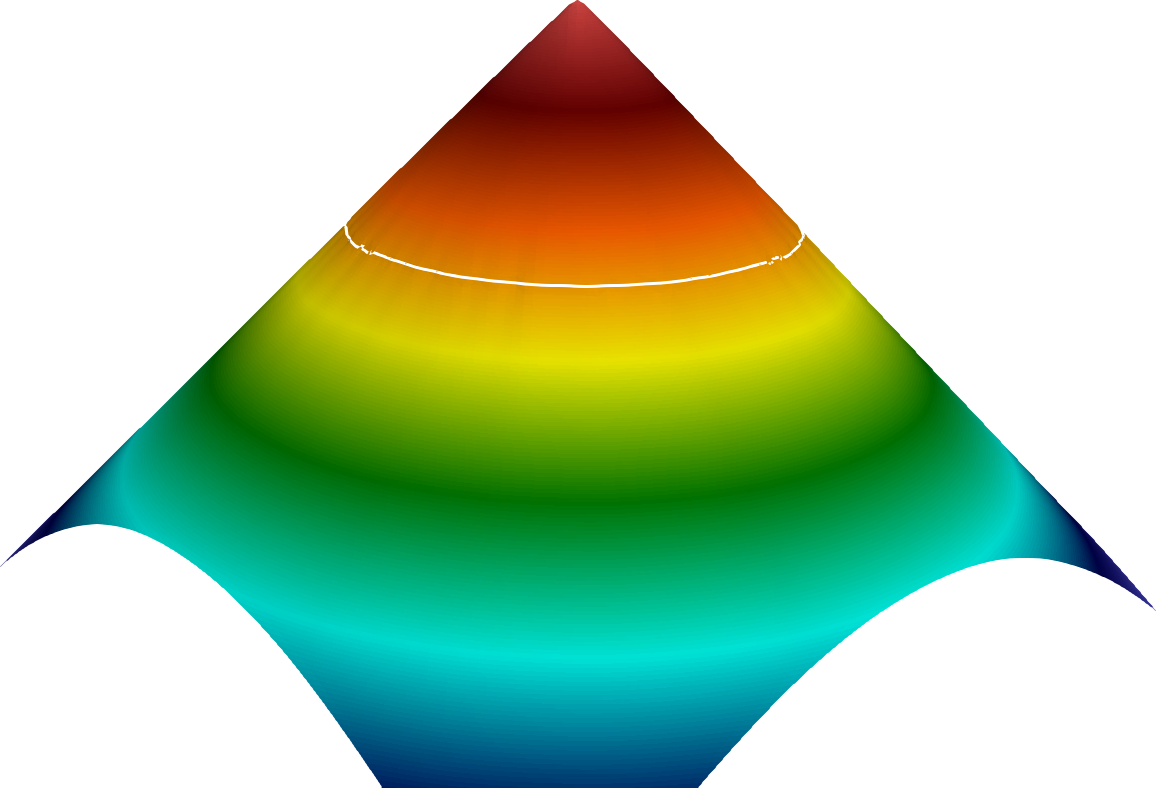}
  \label{fig:disc_ICs_phi0}}

\caption{Redistancing from a discontinuous configuration. We show the (a) initial state,
  (b) the elements containing the interface and (c)-(f) the different steps of this process.
  In all figures we plot the interface in solid white.
  The mesh is uniform with mesh size given by $h=\frac{1}{160}$.}
\end{figure}

The conservative level set model \eqref{level_set_model} requires an initial condition
$\phi(\bfx,t=0)$ to be given by an SDF. From a practical point of view
this might be inconvenient in some situations. Instead, we propose to start the algorithm by considering
an initial configuration given by a discontinuous function $\hat{\phi}\in\{-1,0,1\}$
representing the two phases and re-distance it to obtain an SDF. Due to the fact that the underlying equations for the SDF are nonlinear and the proposed initial conditions are far from the root, we use a robust multi-stage approach to ensure convergence and accuracy. This somewhat complex approach is used only at initialization and requires no user input beyond the initial phase configuration. The main objectives is a practical and robust SDF calculation on unstructured meshes.
See figure \ref{fig:disc_ICs_phiHat} where we consider $\Omega=(0,1)^2$ and
\begin{align}\label{prob_disc_ICs}
  \hat{\phi}=
  \begin{cases}
    1, & \text{ if } ~ \dist(\bfx,\bfx_0) < R, \\
    0, & \text{ if } ~ \dist(\bfx,\bfx_0) = 0, \\
    -1,& \text{ if } ~ \dist(\bfx,\bfx_0) > R, 
  \end{cases}
\end{align}
where $\bfx_0=(0.5,0.5)$ and $R=0.25$.
First we project $\hat{\phi}$ onto the finite element space. Afterwards, 
the main algorithm is based on solving the Eikonal equation several times.
The first time we re-distance the solution away from the interface. Then we concentrate
on the cells containing the interface. And finally we perform a global redistancing
to polish the result. We now provide details about this process. 

{\it Step 1: Lumped $L^2$-projection}.
We start by doing a lumped $L^2$-projection to obtain 
$\hat{\phi}_h=\sum_j \hat{\Phi}_j w_j(\bfx)$ where,
upon defining $m_i=\int_\Omega w_i d\bfx$,
$\hat{\Phi}_i=\frac{h_e}{m_i}\int_\Omega \hat{\phi} w_i d\bfx$.
Here $h_e$ is a characteristic element size;
e.g., $h_e=\frac{1}{2}\left[\max_{\bfx\in\Omega} h(\bfx)+\min_{\bfx\in\Omega}h(\bfx)\right]$.
Note that due to the scaling by $h_e$, $\hat{\phi}_h$ has units of distance.
This step is important to introduce some dissipation into the initial condition.
In figure \ref{fig:disc_ICs_phi_hHat} we show $\hat{\phi}_h$. 

{\it Step 2: Redistancing away from the interface}.
Given $\hat{\phi}_h\in X_h^1$, we solve the (viscous) Eikonal equation away from the interface
to obtain $\phi_h^*\in X_h^1$ as follows:
\begin{align}\label{disc_ICs_step2}
  \int_\Omega \Seps(\hat{\phi}_h)\left[|\nabla\phi^{*}_h|-1\right] w(\bfx) d\bfx
  + c\int_\Omega h(\bfx)\nabla\phi_h^*\cdot\nabla w d\bfx
  + \int_\Omega \tau [\nabla\phi_h^{*}-\bfq_h(\phi_h^{*})]\cdot\nabla w(\bfx)d\bfx = 0,
  & & \forall w(\bfx)\in X_h^1,
\end{align}
During this process, we freeze the DOFs associated with the interface
$\hat{\Gamma}=\{\bfx\in\Omega\st \hat{\phi}_h=0\}$ by imposing strongly
$\Phi_i^*=\hat{\Phi}_i, ~ \forall i\in\I(\hat{\Gamma})$.
The first term in \eqref{disc_ICs_step2} is the consistency term w.r.t. the Eikonal equation.
In the second term $c=\bigO(1)$, we use $c=0.1$, controls the amount of background dissipation. 
We introduce this dissipation to aim the solution to converge to the viscous solution of the Eikonal equation.
The third term acts as non-linear stabilization and penalization from the distance function, see remark \ref{remark:disc_ICs_stab_term}.
In \S \ref{sec:clsvof_normal_reconstruction} we provide details about the
discretization of $\bfq_h(\cdot)$.
In figure \ref{fig:disc_ICs_phi_star} we show the result of this step. 

{\it Step 3: Redistancing close to the interface}.
Now we re-distance the solution close to the interface.
Given $\phi^*_h\in X_h^1$, we solve the Eikonal equation to obtain
$\phi^{**}_h\in X_h^1$ as follows:
\begin{align}\label{disc_ICs_step3}
  \int_\Omega \Seps(\phi^*_h)\left[|\nabla\phi^{**}_h|-1\right] w(\bfx) d\bfx
  + c\int_\Omega h(\bfx)\nabla\phi_h^{**}\cdot\nabla w d\bfx
  +\int_\Omega \tau [\nabla\phi_h^{**}-\bfq_h(\phi_h^{**})]\cdot\nabla w(\bfx)d\bfx = 0,
  & & \forall w(\bfx)\in X_h^1.
\end{align}
During this process, we impose strongly $\Phi_i^{**}=\Phi_i^{*}, ~ \forall i\in\I(\Omega)\setminus\I(\hat{\Gamma})$; 
i.e., we allow changes only on the DOFs associated to the interface.
By doing this, we aim to limit drastic movement of the interface. 
We perform only one Newton iteration of this step.
In figure \ref{fig:disc_ICs_phi_sstar} we plot $\phi_h^{**}$.

{\it Step 4: Global redistancing}. 
Finally we obtain $\phi_h^0\in X_h^1$ by solving 
\begin{align}\label{disc_ICs_step4}
  \begin{split}
  \alpha m_i^{\Gamma^{**}}(\Phi_i^0-\Phi_i^{**})
  +\int_\Omega \Seps(\phi^{**}_h)\left[|\nabla\phi^{0}_h|-1\right] w(\bfx) d\bfx
  +c\int_\Omega h(\bfx)\nabla\phi_h^{0}\cdot\nabla w d\bfx
  \\+\int_\Omega \tau [\nabla\phi_h^{0}-\bfq_h(\phi_h^{0})]\cdot\nabla w(\bfx)d\bfx = 0,
   & \qquad \forall w(\bfx)\in X_h^1.
  \end{split}
\end{align}
The first term, inspired by \cite{basting2014optimal}, is a penalization term to prevent
large changes in the interface. 
Here $\alpha\in\mathbb{R}$ is a penalization constant, we use $\alpha=10^9$,
and $m_i^{\Gamma^{**}}=\int_\Omega \delta_\epsilon(\phi^{**}_h) w_i d\bfx\approx
\int_{\Gamma^{**}} w_i d\bfs$, where $\Gamma^{**}=\{\bfx\in\Omega\st \phi_h^{**}=0\}$.
We perform only one Newton iteration of this step.
We show $\phi_h^0$ in figure \ref{fig:disc_ICs_phi0}. 

  To solve equations \eqref{disc_ICs_step2}-\eqref{disc_ICs_step4}
  we use a quasi-Newton method. Here we show only the
  approximate Jacobian of equation \eqref{disc_ICs_step4}.
  The approximate Jacobian at the $k$-th Newton iteration is given by 
  \begin{align*}
    J_{ij}^k := \alpha m_i^{\Gamma^{**}}\delta_{ij}
    +\int_\Omega \Seps(\phi_h^{**})\frac{(\nabla\phi_h^k\cdot\nabla w_j)}{|\nabla\phi^k|+10^{-10}}w_i d\bfx
    +\int_\Omega \left[ch(\bfx) + \tau\right]\nabla w_i\cdot\nabla w_j d\bfx,
  \end{align*}
  where $\delta_{ij}$ is the Kronecker delta function.
  Note that, for simplicity, we drop the part related to $\bfq_h(\phi_h^0)$.

\begin{remark}[About the stabilization term]\label{remark:disc_ICs_stab_term}
  The last term in equations \eqref{disc_ICs_step2}-\eqref{disc_ICs_step4}
  acts as stabilization. Indeed
$\nabla\phi-\bfq(\phi)\approx \left(1-\frac{1}{|\nabla\phi|}\right)\nabla\phi$.
Therefore, this term behaves like the weak discretization of a Laplacian term
with coefficient given by the residual of the Eikonal equation.
This idea is similar to the high-order stabilization proposed in \cite[\S 5]{lohmann2017flux} and the nonlinear residual-based variational multiscale stabilization used in \cite{kees2011conservative}.
The parameter $\tau$ controls the strength of the stabilization.
We use $\tau=\frac{h(\bfx)}{2}|{\bf V}|$ with $|{\bf V}|=1$ which is reasonable since
${\bf V}=\Seps(\phi)\frac{\nabla\phi}{|\nabla\phi|}$ is the redistancing velocity.
Note that this stabilization term is similar to the regularization and penalization terms in
\eqref{level_set_model_eqn1}. Moreover, the terms in \eqref{level_set_model_eqn1}
can be interpreted as nonlinear stabilization where the coefficient is given not by
the residual of the level set equation but by the residual of the
Eikonal equation. See \cite[Remark 3.2.2]{quezada2018monolithic}.
\end{remark}

\begin{remark}[About the sign of the stabilization term]
  Note that $\left(1-\frac{1}{|\nabla\phi|}\right)<0~$ if $~|\nabla\phi|<1$.
  In this situation anti-diffusion is applied which helps in the convergence to $|\nabla\phi|=1$.
  This is similar to parabolic redistancing by \cite{chan1999nonlinear} and
  elliptic redistancing by \cite{basting2014optimal}.
  Moreover, equations \eqref{disc_ICs_step2}-\eqref{disc_ICs_step4}
  can be seen as a blend between
  hyperbolic and elliptic redistancing. By choosing $\tau=\bigO(h)$ we favor hyperbolic
  redistancing. 
\end{remark}

We consider again $\hat{\phi}$ to be given by \eqref{prob_disc_ICs}
and perform a convergence test. The results are shown in table \ref{table:convergence_disc_ICs}. 

\begin{table}[!ht]
  \centering
  \begin{tabular}{|c|c|c|c|} \hline
    h & N-DOFs & $||\phi_h^0-\phi^\text{exact}||_{L^2(\Omega)}$ & Rate \\ \hline
    2.50E-2 & 1,681    & 6.75E-3 & --   \\ \hline
    1.25E-2 & 6,561    & 2.68E-3 & 1.33 \\ \hline
    6.25E-3 & 25,921   & 1.15E-3 & 1.22 \\ \hline
    3.12E-3 & 103,041  & 5.88E-4 & 0.96 \\ \hline
    1.56E-3 & 410,881  & 3.05E-4 & 0.94 \\ \hline
    7.81E-4 & 1,640,961 & 1.28E-4 & 1.24 \\ \hline
  \end{tabular}
  \caption{Convergence of redistancing from discontinuous configuration.    
  \label{table:convergence_disc_ICs}}
\end{table}

\subsection{Full discretization of the conservative level set method}

\subsubsection{$\mathcal{C}^0$ normal reconstruction}\label{sec:clsvof_normal_reconstruction}
The vector field $\bfq_h(\phi_h)$ in \eqref{level_set_model_eqn2} approximates
a normal field to the interface. We follow \cite{quezada2018monolithic} and consider
a weighted lumped $L^2$-projection given by
\begin{subequations}\label{l2_projection}
\begin{align}
  \bfq_h(\phi_h)=\sum_j \bfQ_j(\phi_h) w_j(\bfx),
\end{align}
where the components of $\bfQ_j\in\mathbb{R}^d, \ j\in\I(\Omega)$ are
calculated using
\begin{align}
  Q_i^{(k)}(\phi_h) = \frac{\int_\Omega\partial_k\phi_h w_i d\bfx}
  {\int_\Omega\sqrt{|\nabla\phi_h|^2+\delta^2}w_i d\bfx}.
\end{align}
\end{subequations}

\subsubsection{Redistancing pre-stage}\label{sec:clsvof_pre_stage}
We begin each time step with a pre-redistancing stage using the algorithm in
\S\ref{sec:disc_initialization}. The motivation behind this step is to introduce
some hyperbolicity into the redistancing process. Our aim is prompt 
the redistancing to emanate from the interface.
In our experience, this is not always
possible with the elliptic redistancing, and thus not always possible with the penalization embedded
in \eqref{level_set_model_eqn1}.
Given the solution $\phi^n_h$ at time $t^n$ we find $\phi^*_h\in X_h^1$ such that 
\begin{align}\label{pre_redist_stage}
  \alpha m_i^{\Gamma^n}(\Phi_i^*-\Phi_i^n)
  +\int_\Omega \Seps(\phi^n_h)\left[|\nabla\phi^*_h|-1\right] w(\bfx) d\bfx
  +\int_\Omega \tau [\nabla\phi^*_h-\bfq_h(\phi^n_h)]\cdot\nabla w(\bfx)d\bfx &= 0,  & & \forall w(\bfx)\in X_h^1,
\end{align}
where $m_i^{\Gamma^n}=\int_\Omega \delta_\epsilon(\phi^n_h)w_i d\bfx$, 
$\alpha=10^9$ and $\tau$ is given as in remark \ref{remark:disc_ICs_stab_term}.
To freeze the interface $\Gamma^n=\{\bfx\in\Omega \st \phi_h^n=0\}$,
we can set $\alpha=0$ and impose strongly
$\Phi_i^*=\Phi_i^n, ~ \forall i\in\I(\Gamma^n)$.
In this case, contrary to \S \ref{sec:disc_initialization},
we keep $\bfq_h$ independent of the solution $\phi^*_h$. This is done to
avoid extra computational effort during this pre-stage.

\subsubsection{Discretization of the level set equation}\label{sec:clsvof_discretization}
In \cite{quezada2018monolithic} it is noted that one can use linear continuous Galerkin
finite elements with no extra stabilization provided that
the advection term is treated implicitly.
We follow this idea and solve model \eqref{level_set_model_eqn1} via

\begin{align}\label{level_set_model_CG}
  \begin{split}
    R(\phi_h^{n+1},w):=\int_\Omega\frac{\Seps(\phi_h^{n+1})-\Seps(\phi^n_h)}{\dt}wd\bfx
    -\frac{1}{2} \int_\Omega \left[\Seps(\phi_h^{n})\bfv^{n}
      + \Seps(\phi_h^{n+1})\bfv^{n+1}\right] \cdot \nabla w d\bfx \qquad\qquad
    \\+\int_\Omega\lambda[\nabla\phi_h^{n+1}-\bfq_h(\phi_h^*)]\cdot \nabla w d\bfx
    +\frac{1}{2}\int_{\partial\Omega^-} \left[\Seps(\phi_h^{n}) \bfv^{n}
      +\Seps(\phi_h^{n+1})\bfv^{n+1}\right] \cdot \bfn w d\bfs
    = 0, \qquad \forall w\in X_h^1,
  \end{split}
\end{align}
where $\phi_h^*$ is the result of the pre-stage given by \eqref{pre_redist_stage}. 
Note that we discretize the consistency terms in \eqref{level_set_model_eqn1}
via a second order Crank-Nicolson discretization
in time and use a first-order, implicit-explicit discretization for the regularization and
penalization terms respectively.
This approach keeps the method simple and efficient (compared to the two stage
method proposed in \cite{quezada2018monolithic}).
Indeed, by doing a linearization of \eqref{level_set_model_eqn1} around the interface,
the regularization and penalization terms are expected to be $\bigO(h^2)$, see
\cite[Remark 3.2.2]{quezada2018monolithic}.
Therefore, one can expect no harm from being permissive with respect to the order of
approximation of such terms. 
We solve \eqref{level_set_model_CG} via Newton's method. The Jacobian corresponding to
$R(\phi^{n+1}_h,w_i)$ at the $k$-th Newton iteration is given by
\begin{align*}
  \frac{\partial R(\phi_h^k,w_i)}{\partial\Phi_j}=:J_{ij}^k =
  \int_\Omega \Seps^\prime(\phi_h^k)
  \left[\frac{1}{\dt}w_iw_j-\frac{1}{2}\nabla w_i\cdot(\bfv^{n+1}w_j)\right] d\bfx+
  \int_\Omega\lambda
  \nabla w_i\cdot \nabla w_j d\bfx\\
  +\frac{1}{2}\int_{\partial\Omega} \Seps^\prime(\phi_h^k)w_iw_j
  \bfv^{n+1}\cdot\bfn d\bfs.
\end{align*}

\begin{figure}[!htb]
  \includegraphics[scale=0.12]{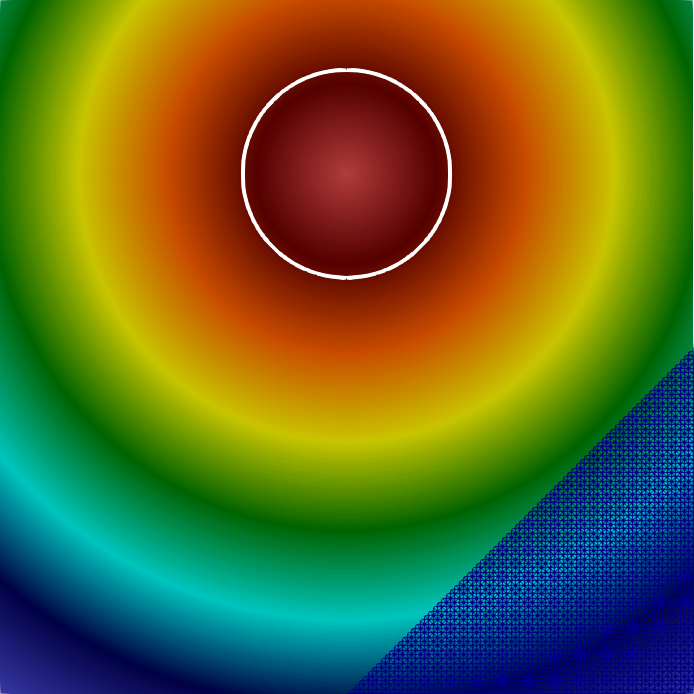}~
  \includegraphics[scale=0.12]{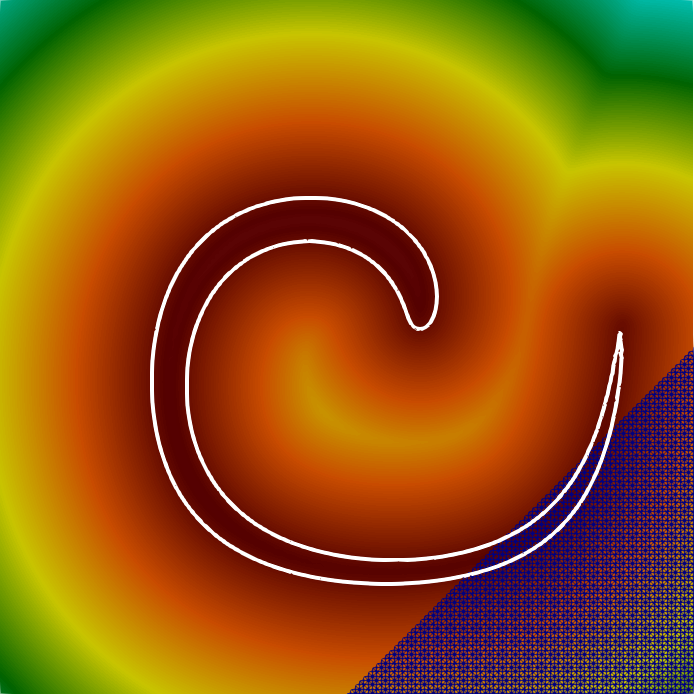}~
  \includegraphics[scale=0.12]{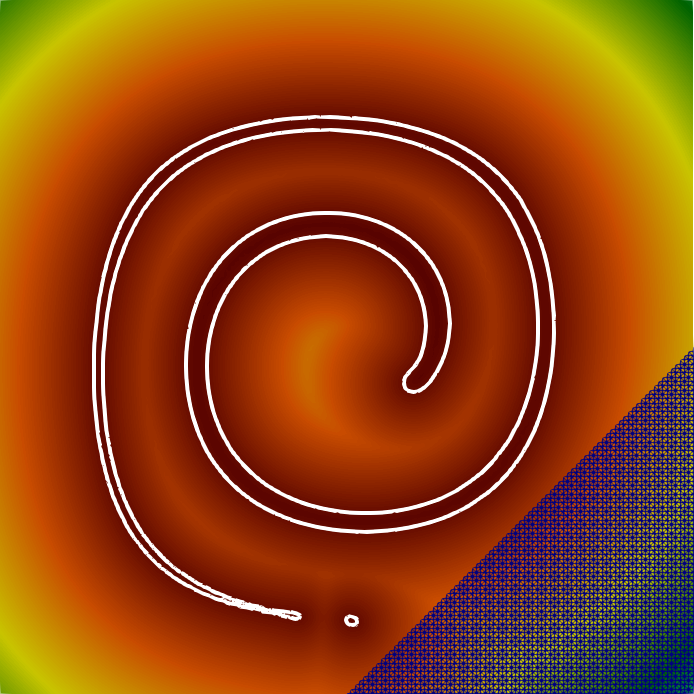}~
  \includegraphics[scale=0.12]{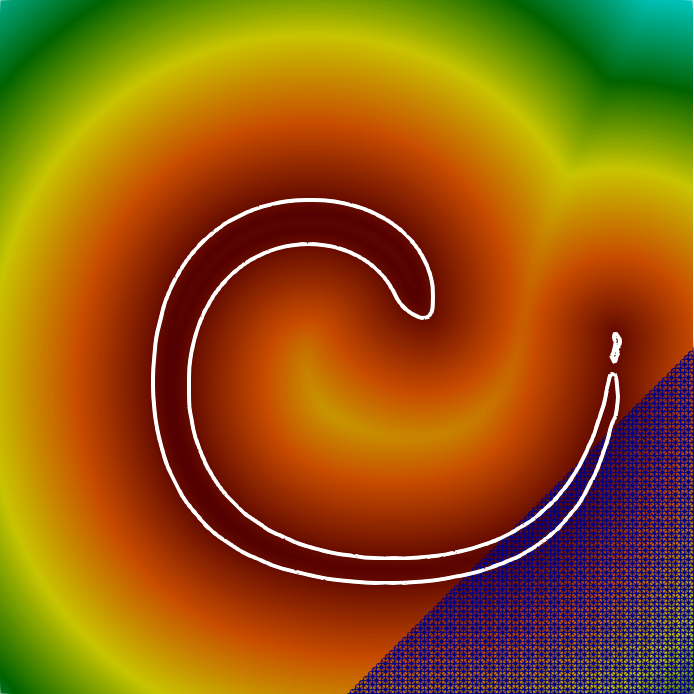}~
  \includegraphics[scale=0.12]{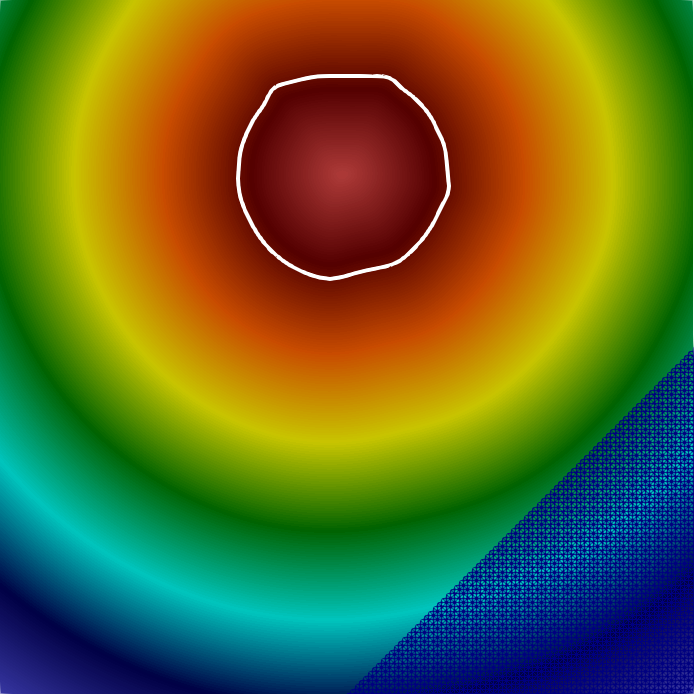}

  \caption{Two dimensional periodic vortex problem.
    We show the solution $\phi$ at $t=0, 2, 4, 6, \text{ and } 8$.
    In all figures we plot the zero contour plot in solid white.
    The mesh, partially displayed in the lower right corner,
    is composed of structured triangles with mesh size given by
    $h=\frac{1}{160}=6.25\times 10^{-3}$.}
\label{fig:clsvof_periodic_vortex}
\end{figure}

We close this section by solving a benchmark in the literature of level sets.
The problem is known as periodic vortex, see \cite{rider1995stretching}.
The domain is given by $\Omega=(0,1)^2$.
The initial condition and velocity field are given by
\begin{subequations}\label{periodic_vortex}
  \begin{align}
    \phi(\bfx,0)&=\pm\dist(\bfx,\Gamma_0),\quad\label{pvortex_init_cond}\\
    \bfv(x,y,t) &=
    \begin{bmatrix}
      -\sin(\pi x)^2\sin(2\pi y)\sin(2\pi t/8) \\
      \sin(2\pi x)\sin(\pi y)^2\sin(2\pi t/8)
    \end{bmatrix},
  \end{align}
\end{subequations}
where $\Gamma_0:=\{(x,y)\in\Omega \st (x-x_c)^2+(y-y_c)^2=r^2\}$
is a circle of radius $r=0.15$ centered at
$(x_c, y_c)=(0.5,0.75)$. We select the positive distance
in \eqref{pvortex_init_cond} if $\bfx=(x,y)$
is inside the circle $\Gamma_0$ and the negative distance otherwise.
In figure \ref{fig:clsvof_periodic_vortex} we show the solution at different
times, and in table \ref{table:convergence_clsvof_periodic_vortex}
we show the errors and convergence rates in the $L^1$- and $L^2$-norms for
different refinements. 
Additionally, we compute some of the metrics in \cite[\S 5]{quezada2018monolithic} given by 
\begin{subequations}\label{metrics}
  \begin{align}
    \text{I}_\text{err}(\phi_h) &= \frac{1}{L}||\Heps(\phi(\bfx))-\Heps\left(\phi_h\left(\bfx,t^{n+1}\right)\right)||_{L^1(\Omega)}, \\
    \text{V}_\text{err}(\phi_h) &= \frac{1}{\int_\Omega H(\phi_h(\bfx,0))d\bfx}
    \left|\int_\Omega\left[H(\phi_h(\bfx,0))-H(\phi_h(\bfx,t))\right]d\bfx\right|, \\
    \text{D}_\text{err}(\phi_h) &= \frac{1}{2}\int_\Omega(|\nabla\phi_h|-1)^2 d\bfx,
  \end{align}
\end{subequations}
where $L$ is the $(d-1)$-dimensional measure of the zero level set $\Gamma(0)$,
and $H(\cdot)$ is the sharp Heaviside function.
The quantities $\text{I}_\text{err}(\phi_h)$ and $\text{V}_\text{err}(\phi_h)$ measure
the extent of interface displacements and area/volume difference with respect to the exact SDF and exact Heaviside function, respectively,
see \cite{enright2002hybrid,olsson2005conservative}, while $\text{D}_\text{err}(\phi_h)$
measures the deviation of $\phi_h$ from a distance function.
We present the results in table \ref{table:2D_vortex_metrics}
for five levels of refinement.

\begin{table}[!ht]
  \centering
  \begin{tabular}{|c|c|c|c|c|c|} \hline
    h & N-DOFs & $||\phi_h-\phi^\text{exact}||_{L^1(\Omega)}$ & Rate
    & $||\phi_h-\phi^\text{exact}||_{L^2(\Omega)}$ & Rate \\ \hline
    2.50E-2 & 1,681   & 5.77E-2 & --   & 7.52E-2 & --   \\ \hline 
    1.25E-2 & 6,561   & 1.87E-2 & 1.62 & 2.10E-2 & 1.84 \\ \hline 
    6.25E-3 & 25,921  & 4.51E-3 & 2.05 & 4.78E-3 & 2.13 \\ \hline
    3.12E-3 & 103,041 & 1.33E-3 & 1.75 & 1.40E-3 & 1.77 \\ \hline
    1.56E-3 & 410,881 & 3.82E-4 & 1.80 & 3.53E-4 & 1.98 \\ \hline         
  \end{tabular}  
  \caption{Convergence of two dimensional periodic vortex problem.
    \label{table:convergence_clsvof_periodic_vortex}}  
\end{table}

\begin{table}[!ht]
  \centering
  \begin{tabular}{|c|c|c|c|c|} \hline
    $h$ & N-DOFs & $\text{I}_\text{err}$ & $\text{V}_\text{err}$ & $\text{D}_\text{err}$  \\ \hline
    2.50E-2 & 1,681   & 6.08E-2 & 2.26E-2 & 2.48E-3 \\\hline
    1.25E-2 & 6,561   & 1.82E-2 & 7.68E-3 & 4.10E-4 \\\hline
    6.25E-3 & 25,921  & 4.51E-3 & 3.76E-3 & 7.12E-5 \\\hline
    3.12E-3 & 103,041 & 1.30E-3 & 1.39E-3 & 1.91E-5 \\\hline
    1.56E-3 & 410,887 & 3.60E-4 & 1.74E-4 & 8.42E-6 \\\hline    
  \end{tabular}
  \caption{
    Error metrics \eqref{metrics} for the two dimensional periodic vortex problem.
    The initial length of the interface is $L\approx 9.42\times 10^{-1}$.
    \label{table:2D_vortex_metrics}}
\end{table}

\section{Incompressible Navier-Stokes solver}\label{sec:navier-stokes}
The incompressible, non-conservative form of the Navier-Stokes equations 
with variable material parameters is given by
\begin{subequations}\label{ns:navier_stokes}
\begin{align}   
    \rho (\partial_t\bfu + (\bfu\cdot\nabla)\bfu)
    		-\nabla\cdot\left(2\mu\varepsilon(\bfu)\right)
    		+\nabla p &= \bff, 
    	& \forall \bfx \in\Omega\\
    \nabla\cdot\bfu & = 0, 
    	& \forall \bfx \in\Omega
\end{align}
\end{subequations}  
where $\varepsilon(\bfu)=\frac{1}{2}\left(\nabla\bfu+\nabla\bfu^T\right)$ 
is the symmetric gradient, 
$\rho$ and $\mu$ are the density 
and viscosity respectively and
$\bfu$, $p$ and $\bff$ are the velocity, pressure and force fields respectively. 
%

\subsection{Time discretization via a projection scheme}\label{sec:projection_scheme}
We consider the second order projection scheme 
for the incompressible Navier-Stokes equations with 
variable density by \cite{guermond2009splitting} and adapt it 
for variable time step sizes.
Let $r := \frac{\Delta t^n}{\Delta t^{n-1}}$. Upon defining 
\begin{subequations} \label{ns:projection_scheme}
\begin{align}		
\BDF(\bfu^{n+1})&:=\frac{(1+2r)\bfu^{n+1}-(1+r)^2\bfu^n + r^2\bfu^{n-1}}{(1+r)\Delta t^n}, \\
\bfu^* &:= (1+r)\bfu^n - r\bfu^{n-1}, \\
p^* &:= (1+r)p^{n}-rp^{n-1},                                                 
\end{align}
the projection method is given as follows:
\begin{align}
  \rho^{n+1}\left[\BDF(\bfu^{n+1})+(\bfu^*\cdot\nabla)\bfu^{n+1}\right]
  -\nabla\cdot(2 \mu^{n+1} \varepsilon(\bfu^{n+1}))
  +\nabla p^* & = \bff^{n+1}, \label{projection_method_momentum_eqn} & \forall & \bfx\in \Omega,\\
  \bfu &= \bfu^\text{BC}, & \forall & \bfx\in\partial\Omega_1, \\
  \varepsilon(\bfu^{n+1})\bfn &= \partial_n\bfu^\text{BC},
  & \forall & \bfx\in\partial\Omega_2, \\
  \bfu\cdot\bfn=0, \qquad [\varepsilon(\bfu^{n+1}) \bfn]_\tau &= 0,
  & \forall & \bfx\in\partial\Omega_3,
\end{align}
where
$\partial\Omega_1\subseteq\partial\Omega$ and
$\partial\Omega_2\subseteq\partial\Omega$
are sections of the boundary where Dirichlet and Neumann boundary
conditions for the velocity field $\bfu$ are applied, respectively,
and $\partial\Omega_3\subseteq\partial\Omega$ corresponds to a slip boundary.
Here $(\cdot)_\tau$ denotes the tangential component of a given vector field. 
The pressure is updated via
\begin{align}
  -\Delta\delta\psi^{n+1} &= -\frac{(1+2r)\min_\bfx (\rho(\bfx,t=0))}
  {(1+r)\Delta t^n}\nabla\cdot\bfu^{n+1},
  & \forall & \bfx\in \Omega, \\
  \delta\psi^{n+1} &= p^\text{BC}-p^n, & \forall & \bfx\in\partial\Omega_4, \\
  \nabla\delta\psi^{n+1}\cdot\bfn|_{\partial\Omega} &= 0,
  & \forall & \bfx\in\partial\Omega\setminus\partial\Omega_4, \\ 
  p^{n+1} &= p^n+\delta\psi^{n+1}, & \forall & \bfx\in\Omega,
  \label{eqn:proj_scheme_pressure_update}
\end{align}
where $\partial\Omega_4\subseteq\partial\Omega$ is the section of the boundary
where the Dirichlet boundary condition $p=p^\text{BC}$ is applied.
We consider the following initial conditions:
\begin{align}
	\bfu(\bfx,t=0) &= \bfu_0(\bfx), \\
    \delta\psi(\bfx,t=0)&= 0, \\
    \quad p(\bfx,t=0) &= p_0(\bfx). 
\end{align}
\end{subequations}
See \S\ref{sec:numerical_examples} for details about common initial and boundary conditions. 

\begin{remark}[About the linearization of the momentum equations]
  Note that we follow \cite{guermond2009splitting} and use $\bfu^*$, a second order extrapolation of the velocity field,
  to linearize the momentum equations. Moreover, we use $\bfu^*$ in the next section to incorporate artificial viscosity
  to stabilize the advective term. By doing this, the equation remains linear with respect to the solution at time $t^{n+1}$.  
\end{remark}

\begin{remark}[About accuracy]
  The expected accuracy of this projection scheme is second order in the $L^2$- and $H^1$-norms.
  We refer the reader to \cite{guermond2009splitting} for a convergence study. 
\end{remark}

\subsection{Artificial viscosity}\label{sec:ns_ev}
In this section, we add artificial viscosity to stabilize the advective term in the momentum equations.
To do this we consider \cite{FCT-book,badia2017monotonicity,barrenechea2016edge,guermond2017secondOrderdij}, 
where artificial dissipative operators are introduced to enforce Discrete Maximum Principles (DMP).
%
%
We do not have a DMP requirement; instead, we are interested in adding enough artificial
viscosity to have a robust and well behaved discretization. This extra viscosity must
vanish as the mesh is refined, must behave properly for unstructured meshes and for all
refinements and should not add significant computational cost to the method. In particular,
we are interested in keeping equation \eqref{projection_method_momentum_eqn} linear.
Finally, we aim to have no tunable parameters and to preserve, to the best of our ability, the accuracy
properties of the underlying method.

\subsubsection*{Componentwise smoothness indicator}
We concentrate first on one component of the velocity vector; e.g., in the $x$-component $u_h$,
whose degrees of freedom are denoted by $U$, see \S\ref{sec:finite_element_discretization}.
The first step is to neglect the force, pressure and viscosity terms in \eqref{projection_method_momentum_eqn};
as a result, we get a hyperbolic system, which is more suitable for the theory developed in the
references above. In the rest of this section we consider backward Euler time stepping; nevertheless,
we apply these results to the second order method in \S \ref{sec:projection_scheme}.
By doing this the full discretization becomes
\begin{align}\label{ns_hyperbolic_part}
  M \left(\frac{U^{n+1}-U^n}{\Delta t}\right)+K(\bfu_h^*)U^{n+1} = 0,
\end{align}
where $M$ is the mass matrix and $K$ is the linearized advective matrix, whose
components are given by $M_{ij}=\int_\Omega w_i(\bfx) w_j(\bfx) d\bfx$ and
$K_{ij}=\int_\Omega \left[\bfu_h^*\cdot\nabla w_j(\bfx)\right]w_i(\bfx)d\bfx$ respectively.
%
Now we introduce an artificial dissipative matrix, see e.g. \cite{FCT-book} and references therein,
whose components are given by
\begin{align}\label{DijLow}
  D_{ij}^\text{Low} = \begin{cases}
    \max\left[K_{ij}(\bfu_h^*), ~0, ~K_{ji}(\bfu_h^*)\right], & \text{ if } i\neq j \\
    -\sum_{k\neq i} D_{ik}^\text{Low}, & \text{ if } i=j
    \end{cases}.
\end{align}
The objective behind this idea is to add sufficient conditions to construct $U^{n+1}$ 
as a convex combination of $U^n$. This, along with using a positive lumped mass matrix,
guarantees the DMP property.
By doing this, however, the accuracy of the method is reduced to first order. To recover the second
order accuracy expected from the projection scheme we use a smoothness based indicator,
see \cite{badia2017monotonicity,barrenechea2016edge,guermond2017secondOrderdij}, given by 
\begin{align}\label{beta}
  \beta_{i}^x := 1-\left[\frac{|\sum_j U_j^* - U_i^*|}{\sum_j |U_j^* - U_i^*|+10^{-10}}\right]^2.
\end{align}
The indicator $0\leq \beta_i^x\leq 1$, which is associated with the $x$-component of the velocity field, 
is (close to) one when the solution is smooth and (close to) zero when the solution is non-smooth and at local extrema.

\subsubsection*{Isotropic artificial dissipation}
The smoothness indicators for the other velocity components are computed similarly. Assume the problem is
three dimensional. Then, we define
\begin{align*}
  \beta_i = \min(\beta_i^x, \beta_i^y, \beta_i^z),
\end{align*}
which is subsequently used to reduce the amount of artificial dissipation from \eqref{DijLow}.
This can be done in different ways. In this work, we follow \cite{badia2017monotonicity} and define
\begin{align*}
  D_{ij} := \begin{cases}
    \max\left[(1-\beta_i) K_{ij}(\bfu_h^*), ~0, ~(1-\beta_j) K_{ji}(\bfu_h^*)\right], & \text{ if } i\neq j \\
    -\sum_{k\neq i} D_{ik}, & \text{ if } i=j
    \end{cases}.  
\end{align*}
Finally, we apply this discrete operator to each component of the velocity and
reincorporate the pressure, force and viscosity terms. 
The full discretization for the $x$-th component is given by
\begin{align*}
  M \left(\frac{U^{n+1}-U^n}{\Delta t}\right)+[K(\bfu_h^*) - D + L] U^{n+1} + F^x = 0,
\end{align*}
and similarly for the other components of the velocity. 
Here $L$ is the stiffness matrix and $F^x$ is a vector that accounts for the pressure and force terms.
Let $f^x$ and $\bfn^x$ denote the $x$-th component of the force field and
the $x$-th component of the normal unit vector to the boundary respectively.
Then, the entries of $L$ and $F^x$ are given by
\begin{align*}
  L_{ij} &= \int_\Omega 2\mu^{n+1}\left[\varepsilon(w_j) \cdot \nabla w_i\right] d\bfx, \\
  F^x_i &= -\int_\Omega p_h^* \partial_x w_i(\bfx) d\bfx + \int_{\partial\Omega} p\bfn^x w_i(\bfs) d\bfs -\int_\Omega f^x w_i(\bfx)d\bfx.
\end{align*}
Note that we integrate by parts for the pressure term.

\begin{remark}[About the DMP property]
  When the smoothness indicator $\beta^x$ is constructed based on $U^{n+1}$,
  the mass matrix is lumped such that $m_i=\int_\Omega w_i d\bfx > 0$ 
  and the corresponding $D_{ij}$ matrix is applied to \eqref{ns_hyperbolic_part},
  the DMP property is guaranteed, see e.g. \cite{badia2017monotonicity}.
  By using $U^*$ instead, the DMP property is likely lost.
  In addition, the viscous, pressure and force terms require extra care. 
  Also note that, to avoid loss of accuracy, we do not lump the mass matrix.
  Nevertheless, as mentioned before, we do not aim to satisfy a DMP.
  Instead, our objective is to add robust and reliable artificial viscosity
  to dissipate small spatial scales not resolved by the computational mesh
  in order to obtain a well-behaved method for any mesh and any refinement level while nevertheless converging rapidly. 
\end{remark}

\begin{remark}[About accuracy]
  It is known that the discrete operator $D_{ij}^\text{Low}$ decreases the order of convergence
  to first order even for smooth solutions, see for instance \cite{FCT-book}.
  Second order convergence is recovered except around local extrema when the
  smoothness indicator $\beta$ is used.
  The degeneracy of accuracy around local extrema occurs since $\beta=0$ at local extrema. 
  To overcome this problem, one can consider an extra smoothness sensor based on second derivatives,
  which is suitable for the velocity space $X_h^2$, see \cite{lohmann2017flux}.
  The objective is to identify smooth local extrema to deactivate the artificial dissipation. 
  In \cite{kuzmin2018hpFEM} this idea is employed for different scalar equations.
  By doing this, the authors demonstrate that the high-order accuracy of the underlying method is recovered.
  We do not explore this idea further in this work. 
\end{remark}

\subsection{Surface tension}\label{sec:surface_tension}
In this section, we summarize the work by \cite{hysing2006new} to incorporate
surface tension effects.
Given $\Gamma\subset\Omega$ to be the interface between the two fluids, the goal is to impose 
$[\bfu]|_{\Gamma}=0$ and $-[-p{\bf I}+2\mu\varepsilon(\bfu)]|_{\Gamma}\cdot\bfn^\Gamma=\sigma\kappa\bfn^\Gamma$, where $\bfn^{\Gamma}$ is the unit normal vector to the interface,
$\sigma$ is the surface tension coefficient and $\kappa$ is the curvature of the interface.
These conditions can be imposed by incorporating
\begin{align}\label{ns:surf_tension_force}
  \bff_\text{st}=-\sigma\kappa\bfn^\Gamma|_{\Gamma}
\end{align}
into the left hand side of \eqref{ns:navier_stokes}.
Let $\underline{\nabla} z(\bfx)=\nabla z(\bfx)-[\bfn^\Gamma\cdot\nabla z(x)]\bfn^\Gamma, ~ \forall \bfx\in\Gamma$
be the tangential gradient with respect to $\Gamma$ and
$\underline{\Delta} z(\bfx)=\underline{\nabla}\cdot\underline{\nabla} z(\bfx), ~ \forall \bfx\in\Gamma$
be the Laplace-Beltrami operator. 
Using the fact that $\kappa\bfn^\Gamma=\underline{\Delta}\bfx|_{\Gamma}$, see for instance
\cite{gilbarg2015elliptic}, the weak form of \eqref{ns:surf_tension_force},
after multiplying $\bff_\text{st}$ by a test function $\bfw(\bfx)$
and integrating by parts, is given by
\begin{align}\label{ns:surf_tension_force_weak}
  \bff_\text{st}=\int_\Gamma \sigma\underline{\nabla}\text{Id}_\Gamma(\bfx)\cdot\underline{\nabla}\bfw(\bfx) dS,
\end{align}
where $\text{Id}_\Gamma(\bfx)$ is the identity map on $\Gamma$.
The boundary integral vanishes for the cases of interest in this work. 
Once the force integral related to the surface tension is defined,
the authors in \cite{hysing2006new} propose to treat \eqref{ns:surf_tension_force_weak}
semi-implicitly based on $\text{Id}_\Gamma^{n+1}\approx\text{Id}_\Gamma^n+\Delta t^{n+1}\bfu^{n+1}$
and to approximate the surface integrals by volume integrals via the regularized delta function of
the level set.
Doing this yields
%
\begin{align}\label{ns:surf_tension_force_semi_implicit}
  \bff_\text{st}^{n+1}\approx
  \int_\Omega \sigma\left[\underline{\nabla}\text{Id}_\Gamma^n\cdot\underline{\nabla}\bfw(\bfx)\right]
  \delta_\epsilon\left(\phi_h^{n}\right)d\bfx
  +\Delta t^{n+1}\int_\Omega \sigma\left[\underline{\nabla}\bfu^{n+1}\cdot\underline{\nabla}\bfw(\bfx)\right]
  \delta_\epsilon\left(\phi_h^{n}\right)d\bfx.
\end{align}
Assuming the acute angle condition is satisfied, see for instance \cite{burman2005stabilized},
it is remarked in \cite{hysing2006new} that since the second term in
\eqref{ns:surf_tension_force_semi_implicit} is positive, it
contributes to the stability properties of the momentum equations. In other words, it adds
viscosity to the velocity at the interface $\Gamma$.
See more details about this discretization in the next section. 

\subsection{Full discretization}\label{ns:sec_full_discretization}
Let $\bfu_h, \bfw\in [X_h^2]^d$ and $p_h, \delta\psi_h, \theta \in X_h^\text{1}$.
The full discretization of \eqref{ns:projection_scheme} is given by 
\begin{subequations} \label{ns:galerkin_discretization}
  \begin{align} \label{ns:galerking_discretization_momentum_eqn}
    \begin{split}
\int_\Omega\rho^{n+1}\left[ 1+\frac{(1+r)\Delta t^n}{1+2r}(\bfu_h^*\cdot\nabla) \right]\bfu_h^{n+1}\cdot \bfw d\bfx
+\frac{(1+r)\Delta t^n}{1+2r}\int_\Omega 2\mu^{n+1}\varepsilon(\bfu_h^{n+1}):\varepsilon(\bfw) d\bfx ~ + ~\text{stab.} \\
+\frac{(1+r)\Delta t^n}{1+2r}\Delta t^{n+1}\int_\Omega\sigma
(\underline{\nabla}\bfu_h^{n+1} : \underline{\nabla} \bfw) \delta_\epsilon(\phi_h^n) d\bfx
=\bfb^{n+1}(\bfw),
\end{split}
\end{align}
  where `$\text{stab.}$' refers to the artificial viscosity from \S\ref{sec:ns_ev} and $\bfb^{n+1}$ is given by
  \begin{align}
    \begin{split}
      \frac{1+2r}{(1+r)\Delta t^n}\bfb^{n+1}(\bfw) = \int_\Omega\left[ \rho^{n+1}
    	\left((1+r)^2 \bfu_h^n
    	-r^2 \bfu_h^{n-1}\right) +\bff^{n+1} \right]\cdot\bfw d\bfx
      + \int_\Omega p_h^* \text{div}(\bfw)d\bfx
      \\
      - \int_{\partial\Omega} p_h^*\bfn \cdot \bfw d\bfs
      -\sigma\int_\Omega(\underline{\nabla}\text{Id}_\Gamma^n : \underline{\nabla}\bfw)
      \delta_\epsilon(\phi_h^n) d\bfx 
      + \int_{\partial\Omega_2} 2\mu^{n+1}\partial_\bfn\bfu^\text{BC}\cdot\bfw d\bfs.       
    \end{split}
  \end{align}
Upon defining $a:=\frac{(1+2r)\min_\bfx(\rho(\bfx,t=0))}{(1+r)\Delta t}$, 
the full discretization of the pressure update is given by
\begin{align}
  \int_\Omega \nabla\delta\psi_h^{n+1}\cdot\nabla\theta d\bfx  
  &= a\left(\int_\Omega \bfu_h^{n+1}\cdot\nabla \theta d\bfx
  - \int_{\partial\Omega}(\bfu_h^{n+1}\cdot \bfn) \theta d\bfs\right)
	\label{eqn:proj_scheme_full_disc_pressure_increment}\\
  p_h^{n+1} &= p_h^n + \delta\psi_h^{n+1}.\label{eqn:proj_scheme_full_disc_pressure_update}
\end{align}
\end{subequations}
\begin{remark}[Velocity correction]\label{ns_remark:velocity_correction}
  After the pressure update, it is possible to correct the velocity field 
  to obtain a weakly divergence free velocity field. To do this we can rearrange the 
  terms in \eqref{eqn:proj_scheme_full_disc_pressure_increment} to obtain
  \begin{align*}
    \int_\Omega \left(\bfu_h^{n+1}-\frac{1}{a}\nabla\delta\psi_h^{n+1}\right)\cdot\nabla \theta d\bfx 
    - \int_{\partial\Omega}
    \bfu_h^{n+1}
    \cdot \bfn \theta d\bfs=0,
  \end{align*}
  and redefine the velocity such that
  $\bfu_h^{n+1} \rightarrow \bfu_h^{n+1}-\frac{1}{a}\nabla\delta\psi_h^{n+1}$ in the interior of the domain.
  By doing this it is clear that 
  \begin{align}
    \int_\Omega \bfu_h^{n+1} \cdot\nabla \theta d\bfx 
    - \int_{\partial\Omega} (\bfu_h^{n+1}\cdot \bfn) \theta d\bfs = 0
    \implies
    -\int_\Omega\left(\nabla\cdot\bfu_h^{n+1}\right)\theta d\bfx = 0.
  \end{align}
%
  It is remarked in \cite[\S 3.5]{guermond2006overview}, and references therein,  
  that not doing this correction does not affect the accuracy properties of the method. 
  However, as explained in \S\ref{sec:level_set_method}, 
  the level set model assumes the velocity to be (at least weakly) divergence free. 
  Therefore, we correct the velocity field. 
\end{remark}

\section{Numerical experiments}\label{sec:numerical_examples}
A total of eight test problems were selected to demonstrate the behavior of the method we propose. 
The results were compared qualitatively and quantitatively
versus other results in the literature and versus experimental measurements.
These two-phase flow experiments represent a wide range of applications
ranging from high-viscosity buckling fluids to free surface flows around obstacles.
Our main objective is to test the robustness of the method for a wide range of physical parameters,
surface tension, different type of boundary conditions, external forces and arbitrary mesh refinements.
In particular, we set the numerical parameters to be the same for all problems; in our opinion,
this demonstrates the robustness of the numerical method and suitability of the computational model for engineering applications.

\subsubsection*{Parameters}
The physical parameters are part of the definition of the problem. They are given by the density ($\rho$)
and viscosity ($\mu$) of each phase, the magnitude of the gravity ($g$) and the surface tension coefficient ($\sigma$).
Unless otherwise stated, we consider the metric system and set
\begin{subequations}
\begin{align}\label{material_parameters}
  \rho_W=998.2 ~\left(\frac{\text{kg}}{\text{m}^3}\right), \qquad
  \rho_A=1.205 ~\left(\frac{\text{kg}}{\text{m}^3}\right), \qquad
  \mu_W/\rho_W=1.004\times 10^{-6} ~\left(\frac{\text{m}^2}{\text{s}}\right), \\
  \mu_A/\rho_A=1.5\times 10^{-5} ~\left(\frac{\text{m}^2}{\text{s}}\right), \qquad
  g=9.8 ~\left(\frac{\text{m}}{\text{s}^2}\right), \qquad
  \sigma=72.8\times 10^{-3} ~\left(\frac{\text{kg}}{\text{s}^2}\right),
\end{align}
\end{subequations}
where $W$ and $A$ refer to water and air, respectively. In the rest of this work we omit the units. 
The parameter $\tilde{\lambda}$ in the conservative level set method is the main numerical parameter, see remark \ref{remark:about_lambda}.
In all simulations we use $\tilde{\lambda}=10$.
A less important parameter is $\epsilon$, which is used to regularize Heaviside and Dirac delta functions, see \S\ref{sec:level_set_method}.
During the redistancing process we set $\epsilon=\frac{1}{3}$ and elsewhere we use $\epsilon=\frac{3}{2}$.
%
%
In all problems we let $h_e$ denote a characteristic element size (of the unstructured mesh) and report the number of elements.

\subsubsection*{Commonly used initial and boundary conditions}
In different problems we start `at rest'; i.e., the initial velocity and pressure fields are set to zero.
Likewise, we use similar and common boundary conditions in different problems.
We refer to them as slip, non-slip, open top and inflow boundary conditions.
These common boundary conditions for the velocity are shown in table \ref{table:boundary_conditions},
where $(\cdot)_\tau$ denotes the tangential component of a given vector field. 
The level set does not require boundary conditions when the boundary is set to slip or non-slip.
When the boundary is open and $\bfu\cdot\bfn<0$ we set $\Seps(\phi)=1$; i.e., we let only air through the open boundary.
Finally, when the boundary is of inflow type we set either $\Seps(\phi)=0$ or $\Seps(\phi)=1$ to let
water or air respectively into the domain.

\begin{table}[!h]
  \centering
  \begin{tabular}{|c|c|c|c|}\hline
    Slip & Non-slip & Open top & Inflow \\ \hline
    $\bfu\cdot\bfn=0, \quad [\varepsilon(\bfu)\cdot\bfn]_\tau=0$
    & $\bfu=0$
    & $\varepsilon(\bfu)\cdot\bfn=0, \quad p=0$
    & $\bfu=\bfu^\text{inlet}~$, s.t. $\bfu\cdot\bfn<0$ \\ \hline
  \end{tabular}
  \caption{Commonly used boundary conditions for the velocity. \label{table:boundary_conditions}}
\end{table}

\subsubsection*{Computational framework}
This work was developed using {\sc Proteus} \href{https://proteustoolkit.org}{ (https://proteustoolkit.org)},
a toolkit for computational methods and simulation released under the MIT open source license with source available at \url{https://github.com/erdc/proteus}.
The numerical linear algebra is handled by PETSc, see \cite{petsc-web-page,petsc-user-ref,petsc-efficient}, through petsc4py \cite{dalcin2011parallel}.
Inside {\sc Proteus} 1.5.1, we created a set of files to facilitate the definition of the different problems.
This framework is not part of {\sc Proteus}; nevertheless, to facilitate its use, we provide
a release with this framework incorporated,
see release 1.5.1-mp-r1 (\url{https://github.com/erdc/proteus/releases/tag/1.5.1-mp-r1}).
Our aim is to provide a solid, robust and easy to use open-source computational framework for users interested in two-phase flows.
We do not assume the users have any knowledge on finite elements nor expect the need to adjust the numerical parameters, except for potentially $\tilde{\lambda}$. 
Instead, we expect the user to provide information only about the problem; in particular, the domain, initial and boundary conditions and physical parameters.
We encourage the interested reader to install the {\sc Proteus} release and
run different problems.


\subsection{Two-dimensional experiments}

\subsubsection{Rising bubbles}\label{sec:rising_bubbles}
\begin{figure}[!h]
    \centering
    \subfloat[Test case 1.]
    {
      \includegraphics[scale=0.18]{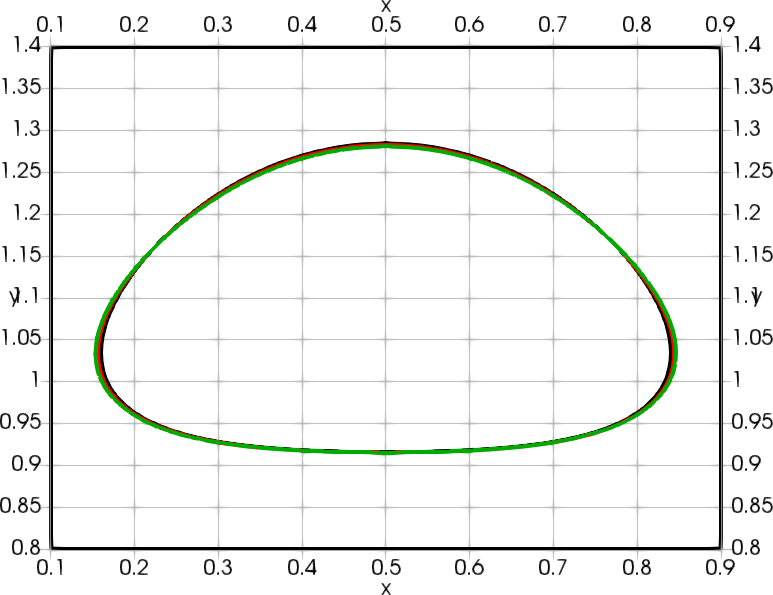}\qquad
      \includegraphics[scale=0.28]{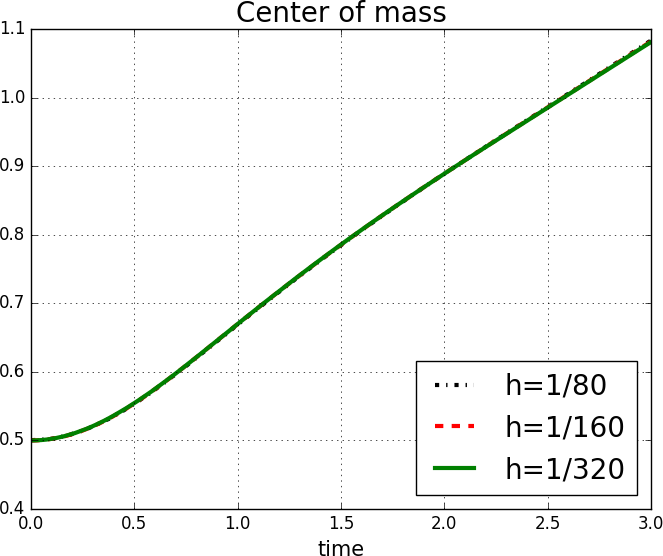}\qquad
      \includegraphics[scale=0.28]{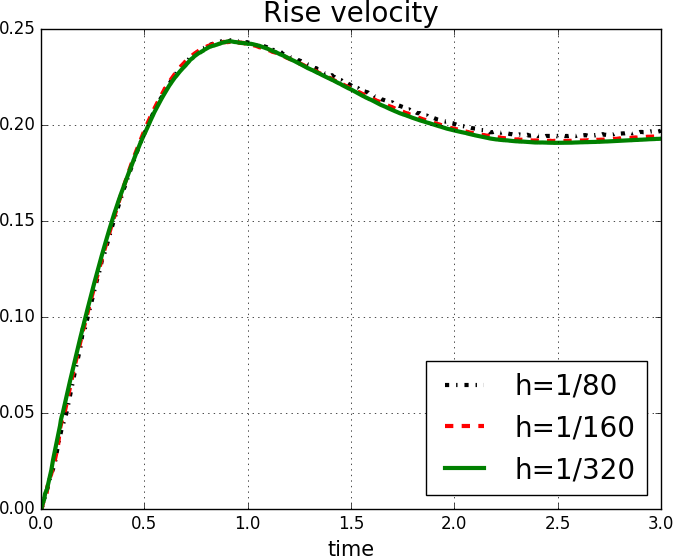}    
    }
    
    \subfloat[Test case 2.]
    {
      \includegraphics[scale=0.18]{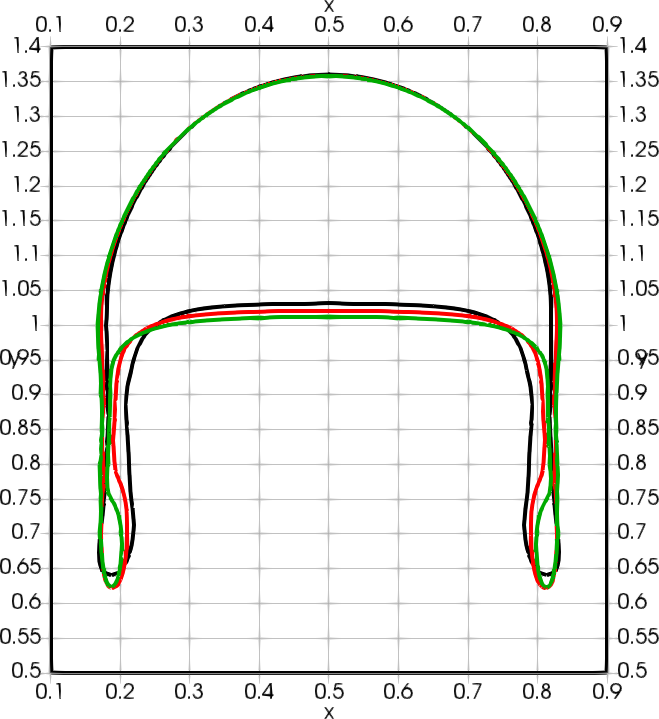}\qquad
      \includegraphics[scale=0.28]{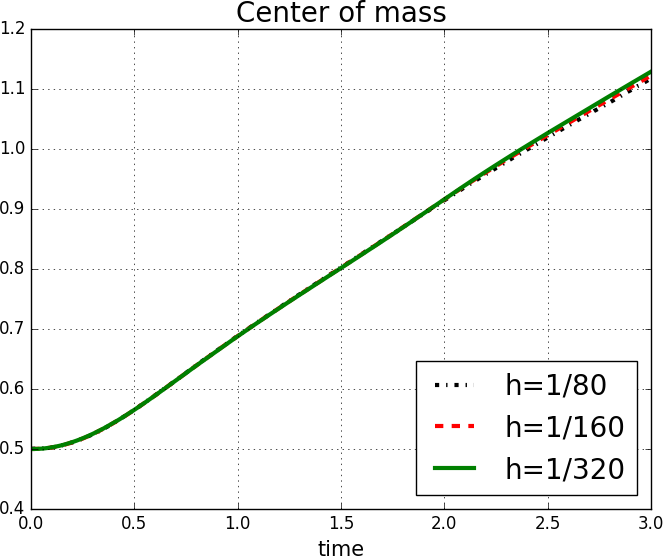}\qquad
      \includegraphics[scale=0.28]{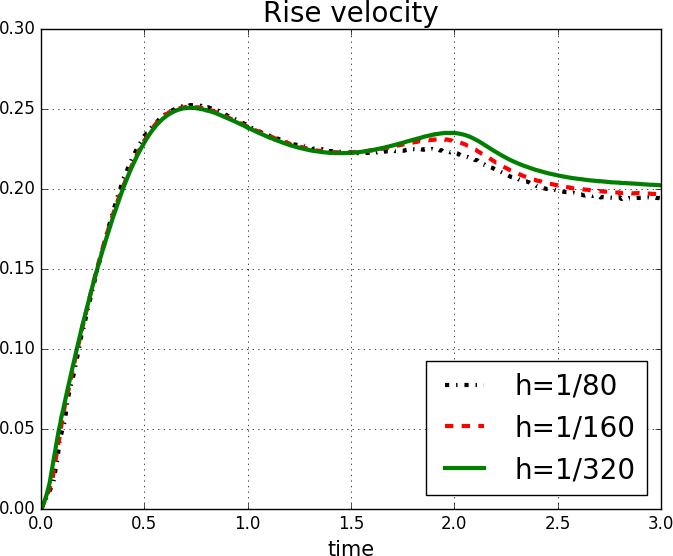}    
    }    
    \caption{
      Rising bubble problems in \cite{hysing2009quantitative}.
      We show (in the left panel) the zero contour plot of the level set at $t=3$,
      (in the middle panel) the evolution of the center of mass
      and (in the right panel) the evolution of the $y$-th component of the velocity.
      We consider three refinement levels with
      (black) $h_e=1/80$, (red) $h_e=1/160$ and (green) $h_e=1/320$. 
      \label{fig:rising_bubble}
    }
\end{figure}

We consider first a bubble of a light fluid at rest immersed in a heavier fluid. 
Due to the action of gravity, the bubble rises through the heavier fluid. 
In this problem, the effect of surface tension is critical to maintain the 
correct shape of the bubble and to obtain the correct position and rising speed. 
We consider the two common test cases in \cite{hysing2009quantitative}. 
The domain of interest is $\Omega=(0,1)\times(0,2)$. We impose non-slip boundary conditions 
on the bottom and top boundaries and slip boundary conditions on the left and right boundaries. 
The initial condition for level set is initialized as explained in \S\ref{sec:disc_initialization}
considering an interface given by $\Gamma=\{(x,t)\in\Omega\st r^2=(x-0.5)^2+(y-0.5)^2, r=0.25\}$.
The material parameters for the different test cases are shown in table
\ref{table:rising_bubble_parameters}.
We consider structured meshes with refinement levels
$h_e=\frac{1}{80}, \frac{1}{160}$ and $h_e=\frac{1}{320}$,
which corresponds to 6480, 25760 and 102720 elements respectively. 
The zero contour plot of the level set and the evolution of 
the center of mass $c_y=\frac{1}{|\Omega_B|}\int_{\Omega_B} y d\bfx$
and rising velocity $r_v=\frac{1}{|\Omega_B|}\int_{\Omega_B} \bfu_y d\bfx$
are shown in figure \ref{fig:rising_bubble}. 
Here $\Omega_B=\{(x,y)\in\Omega\st \phi(x,t)\geq 0\}$ denotes the domain occupied by the bubble. 

\begin{table}[!h]
\centering
	\begin{tabular}{|c|c|c|c|c|c|c|} \hline
        Test case & $\rho_W$ & $\rho_A$ & $\mu_W$ & $\mu_A$ & $g$  & $\sigma$  \\ \hline
        1         & 1000     & 100      & 10      & 1       & 0.98 & 24.5 \\ \hline
        2         & 1000     & 1        & 10      & 0.1     & 0.98 & 1.96 \\ \hline   
    \end{tabular}
\caption{Material parameters for the rising bubble problems. 
			We consider two common test cases. 
		\label{table:rising_bubble_parameters}
		}
\end{table}

\subsubsection{Dam break problem}
\begin{figure}[ht]
  \centering
  \subfloat[Pressure gauge $P_1$.]{
    \begin{tabular}{ccc} 
      \includegraphics[scale=0.185]{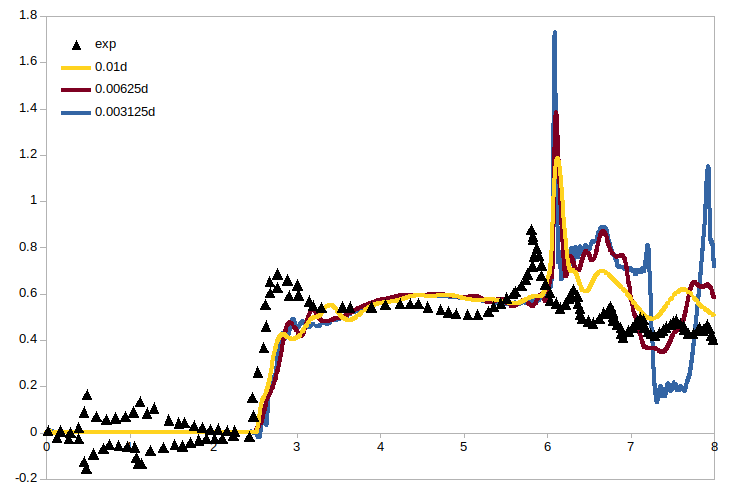}
    \end{tabular}
  }    
  \quad
  \subfloat[Water height gauges $H_1$ and $H_2$.]{
    \begin{tabular}{ccc}
      \includegraphics[scale=0.21]{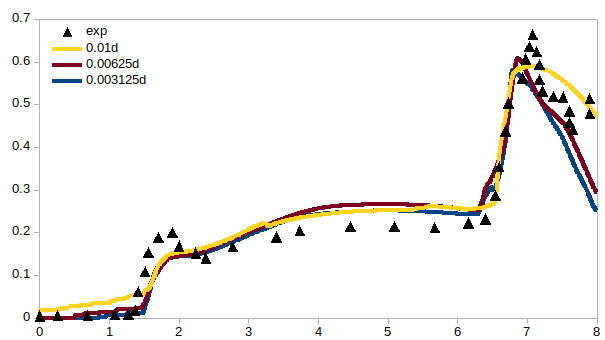} &
      \includegraphics[scale=0.21]{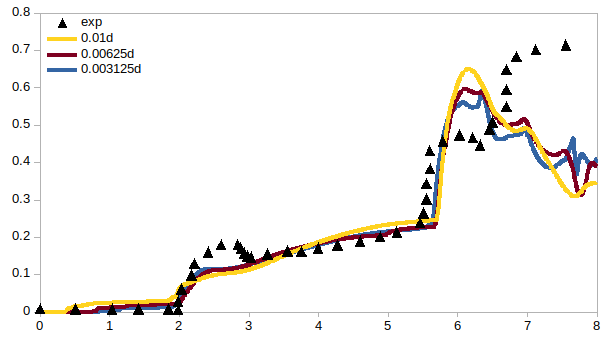}
    \end{tabular}
  }
  \caption{Pressure and water height gauges for the two-dimensional dam break problem.}
  \label{fig:damBreak2D_gauges}
\end{figure}

\begin{figure}[ht]
  \centering
  \includegraphics[scale=0.16]{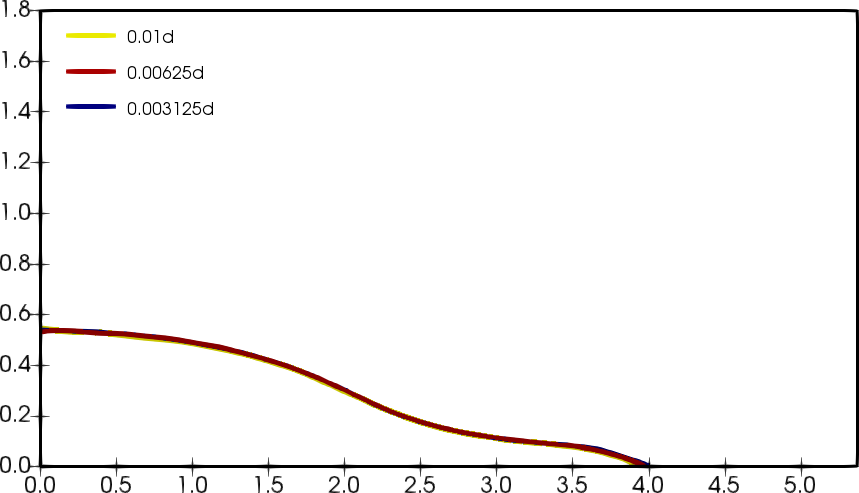}\quad
  \includegraphics[scale=0.16]{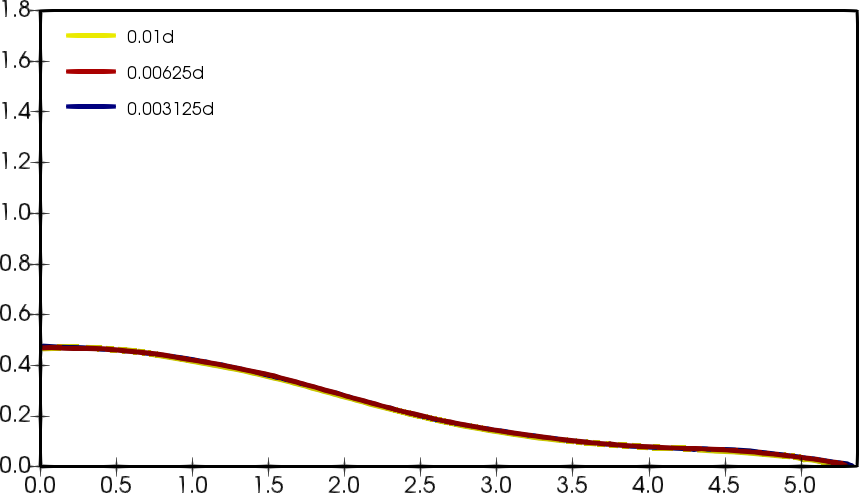}\quad
  \includegraphics[scale=0.16]{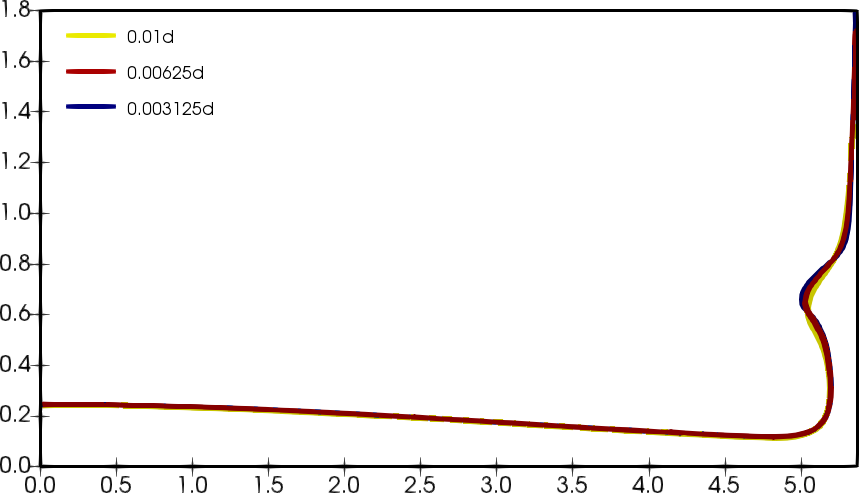}
  
  \includegraphics[scale=0.16]{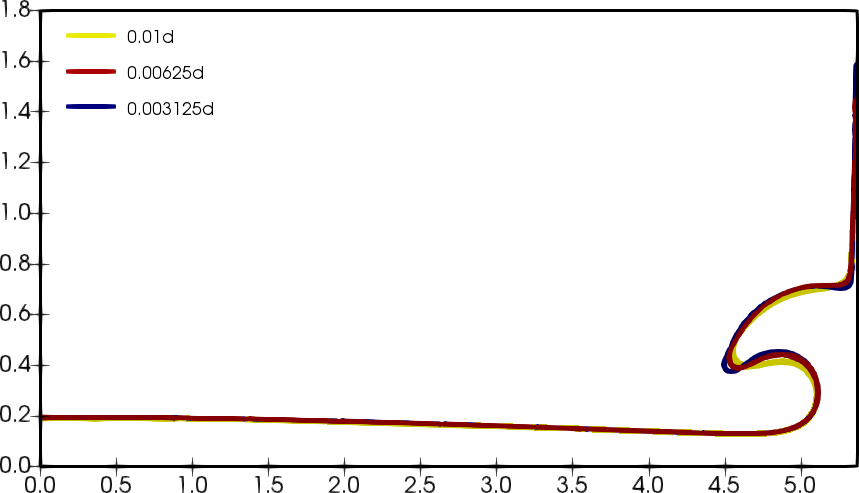}\quad
  \includegraphics[scale=0.16]{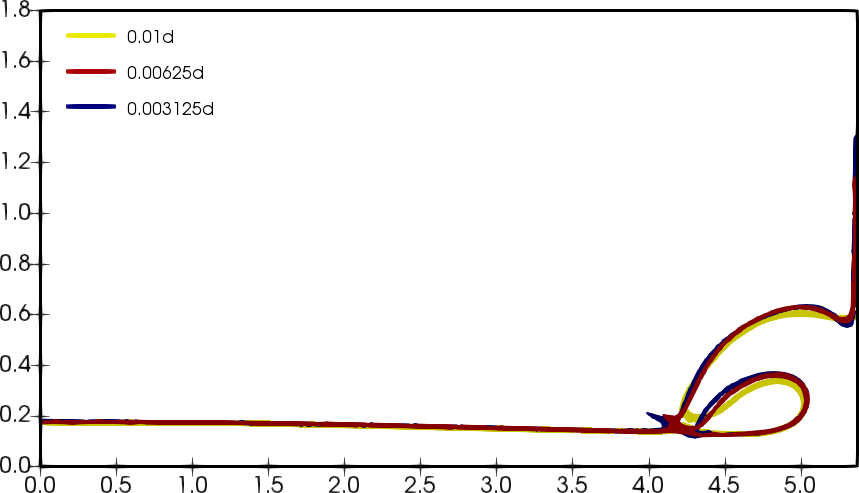}\quad
  \includegraphics[scale=0.16]{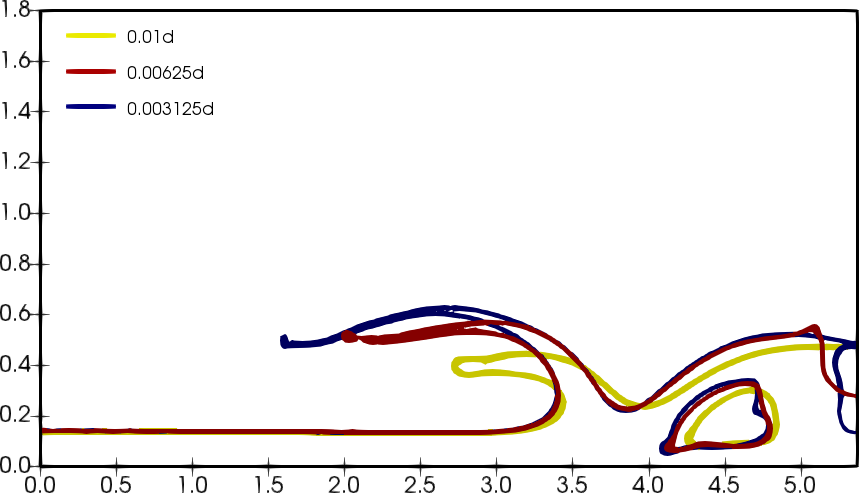}
  \caption{Water-air interface at $\tilde{t}=1.66, ~2.4, ~4.81, ~5.72, ~6.17$ and $\tilde{t}=7.37$.
    Following \cite[figure 17]{colagrossi2003numerical} and to facilitate comparisons, we scale the
    $x$-axis as $\tilde{x}=\frac{x}{H}$ and the time as $\tilde{t}=t\sqrt{g/H}$, where $H=0.6$ is the initial height of the column of water.}
  \label{fig:damBreak2D}
\end{figure}

Here we consider a common two dimensional dam break problem.
In \cite{martin1952part}, the authors considered a similar problem and performed a series of experiments
placing pressure and water height gauges at different locations.
Based on these results, different numerical studies have been performed to validate numerical methods and codes.
In this work, we consider the setup by \cite{colagrossi2003numerical} and reproduce some of their results.
The domain of interest is $\Omega=(0,3.22)\times(0,1.8)$. An initial column of water at rest is located in
$W=\{\bfx\in\Omega \st x\leq 1.2, ~ y\leq 0.6\}$. The rest of the domain ($A=\Omega\setminus W$) is filled with air. 
At time $t=0$ the column of water starts to fall because of the action of gravity.
For this problem, we follow \cite{colagrossi2003numerical} and consider non-slip boundary conditions
everywhere except for the top boundary, which is left open. 
We report the results of only one pressure gauge located at $P_1 = (3.22,0.12)$ and 
two water height gauges located at $H_1=(2.228,y)$ and $H_2=(2.724,y), ~\forall y\in[0,1.8]$. 
Let $d=3.22$ denote the length of the domain.
We consider three unstructured meshes with $h_e=0.01d, ~0.00625d$ and $h_e=0.003125d$, 
which correspond to 17752, 45188 and 181033 triangular elements respectively.
In figure \ref{fig:damBreak2D_gauges} we compare the numerical results (for the different refinements) versus
the experimental measurements for all the gauges.
Finally, in figure \ref{fig:damBreak2D} we consider the same times as
\cite[figure 17]{colagrossi2003numerical} and show the water-air interface.
Note that only in this figure, to facilitate comparisons, we follow the reference and scale the $x$-axis as
$\tilde{x}=\frac{x}{H}$ and the time as $\tilde{t}=t\sqrt{g/H}$, where $H=0.6$ is the initial height of the column of water.

\subsubsection{Buckling flow}
\begin{figure}[ht]
  \centering
  \includegraphics[scale=0.115]{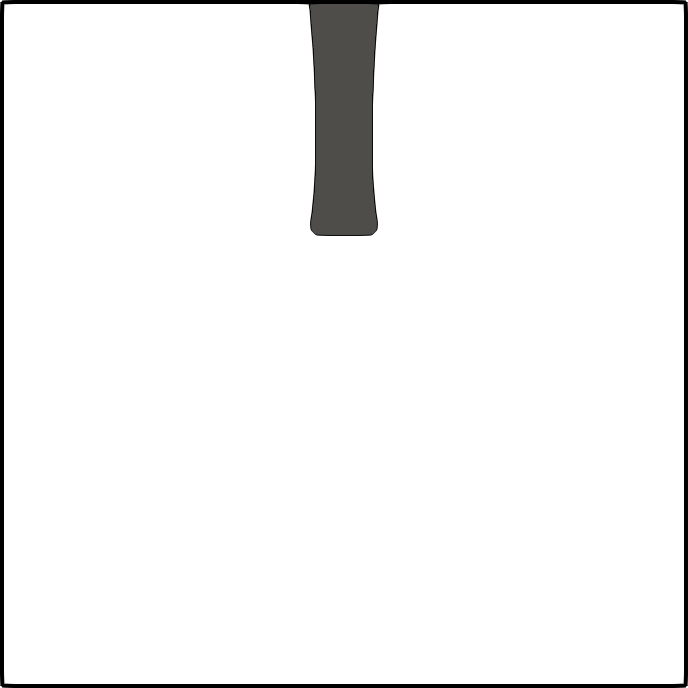}\quad
  \includegraphics[scale=0.115]{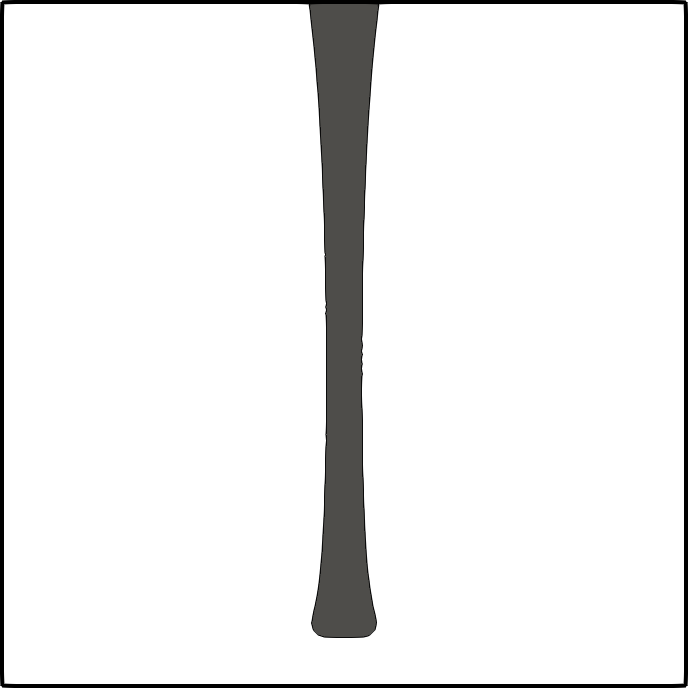}\quad
  \includegraphics[scale=0.115]{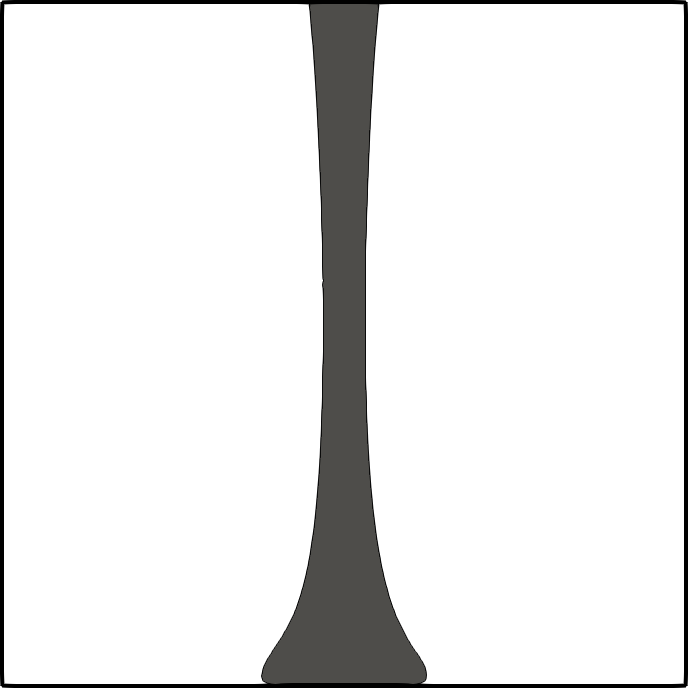}\quad
  \includegraphics[scale=0.115]{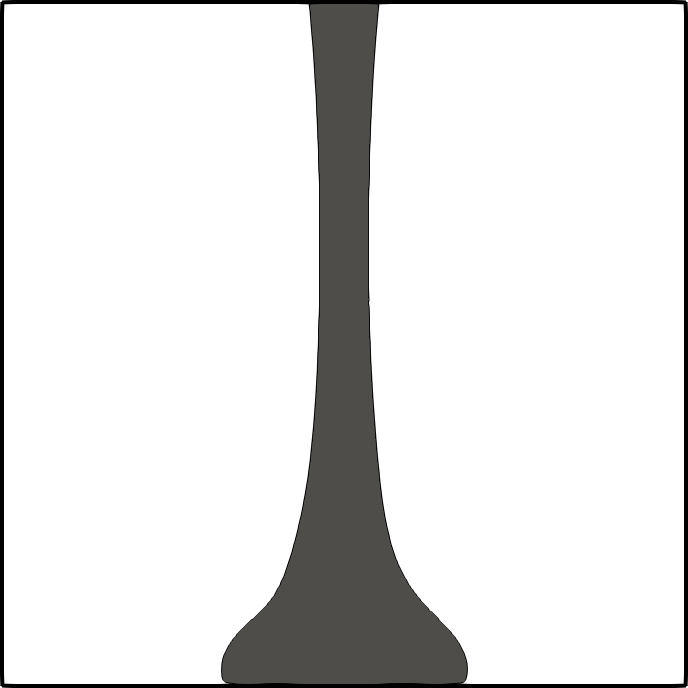}\quad
  \includegraphics[scale=0.115]{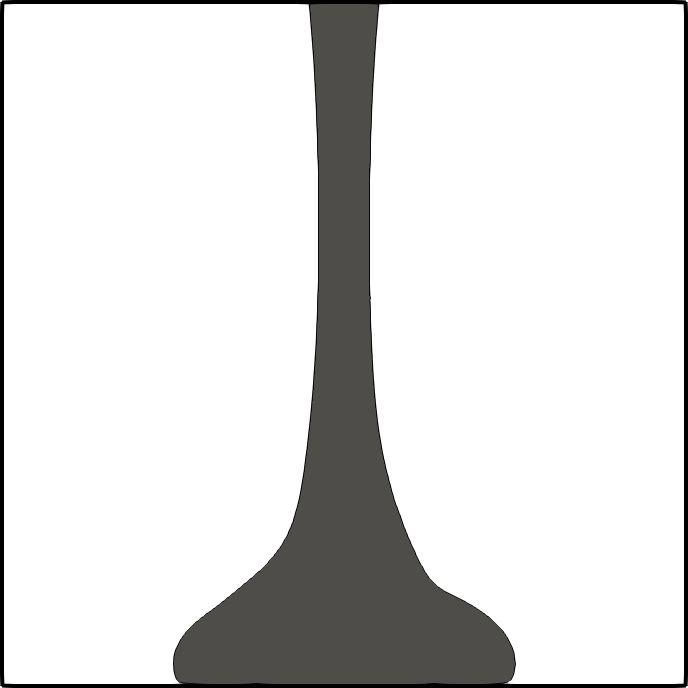}

  ~
  
  \includegraphics[scale=0.115]{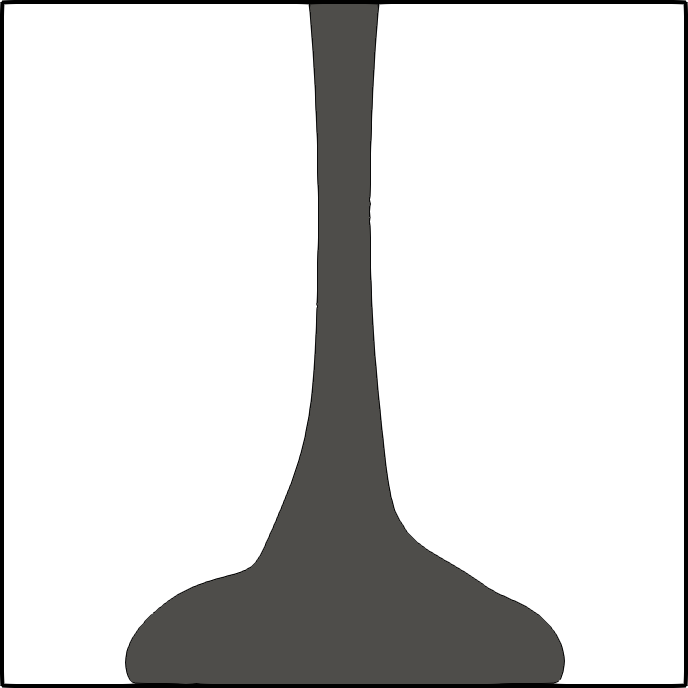}\quad
  \includegraphics[scale=0.115]{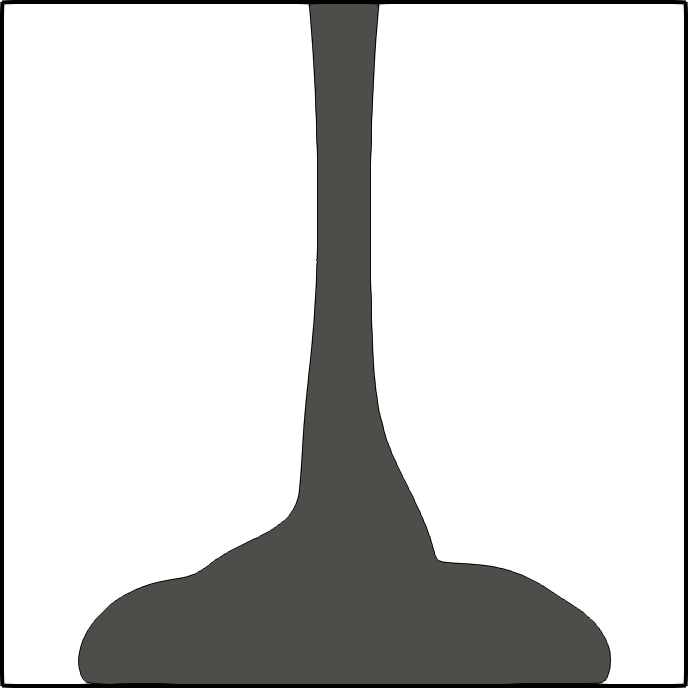}\quad
  \includegraphics[scale=0.115]{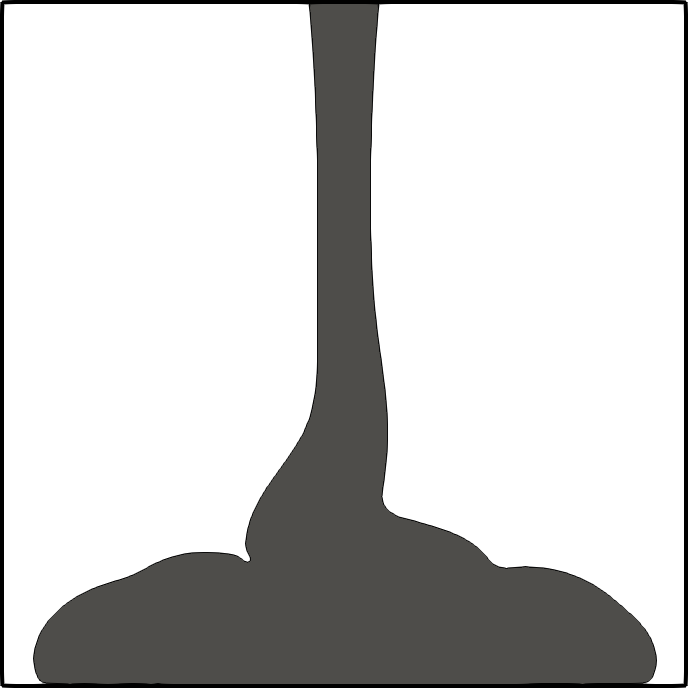}\quad
  \includegraphics[scale=0.115]{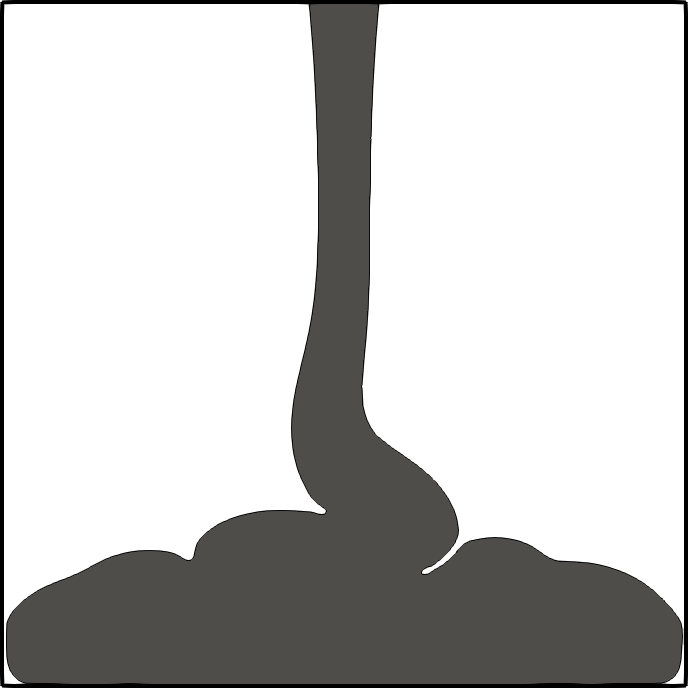}\quad
  \includegraphics[scale=0.115]{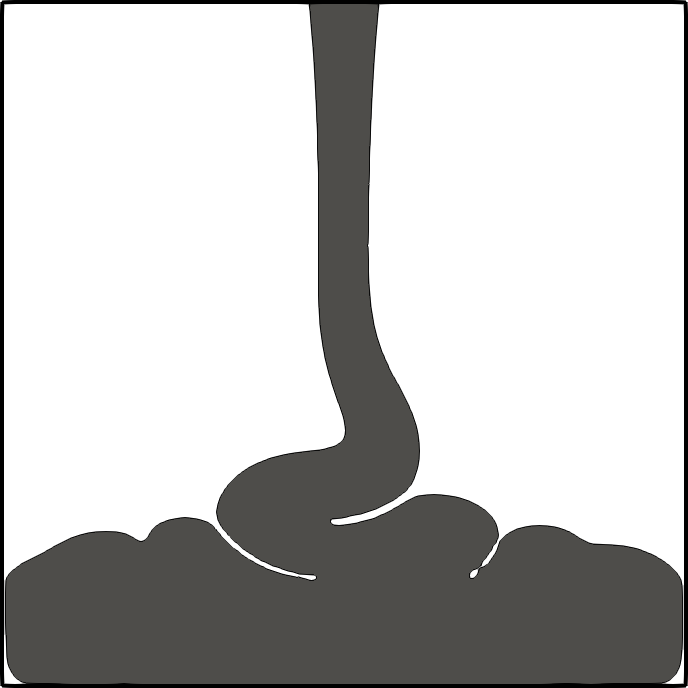}
  \caption{Two-dimensional buckling flow. From top to bottom and left to right,
    we show the heavy fluid phase at
    $t=0.3, ~0.6, ~0.9, ~1.2, ~1.5, ~1.8, ~2.1, ~2.4, ~2.7$ and $~3$.}
  \label{fig:buckling}
\end{figure}

The fluid buckling phenomenon occurs in situations involving thin streams of highly viscous
flow that encounter a boundary or plate. Depending on the Reynold's number and the
cross sectional geometry of the viscous jet, the fluid may fold, coil, or oscillate
antiperiodically. Here, the method's ability to simulate buckling
flow is evaluated against results from previous literature, see e.g.,
\cite{bonito2006numerical, ville2011convected, tome1999numerical, bonito2016numerical}.
The domain $\Omega=(0,1)^2$ is initially occupied by the air phase, which is at rest.
A small inlet on the top boundary of the box, defined by
$I=\{\bfx \in \bar{\Omega} \st |x-0.5| <0.05\}$, is the inflow boundary,
where the velocity is strongly set to $\bfu=(0,-1)$.
The rest of the top boundary is left open.
The non-slip boundary condition is applied in the bottom, left and right boundaries.
The material properties are shown in table \ref{table:buckling_flow}. 
We use an unstructured mesh with a maximum element size of $h_e=5\times10^{-3}$,
which correspond to 126,646 triangular elements.
In figure \ref{fig:buckling} we plot the heavy fluid phase for different times. 

\begin{table}[!h]
\centering
\begin{tabular}{|c|c|c|c|c|c|} \hline
$\rho_W$  & $\rho_A$ & $\mu_W$ & $\mu_A$ & $g$ & $\sigma$  \\ \hline
1800  & 1 & 500 & $2\times 10^{-5}$ & 9.8 & 0 \\ \hline
\end{tabular}
\caption{Material parameters for the buckling flow problem.}
  \label{table:buckling_flow}
\end{table}

\subsubsection{Filling tank}
\begin{figure}[ht]
  \centering
  \subfloat[Test case 1]{
    \begin{tabular}{cccc} 
      \includegraphics[scale=0.13]{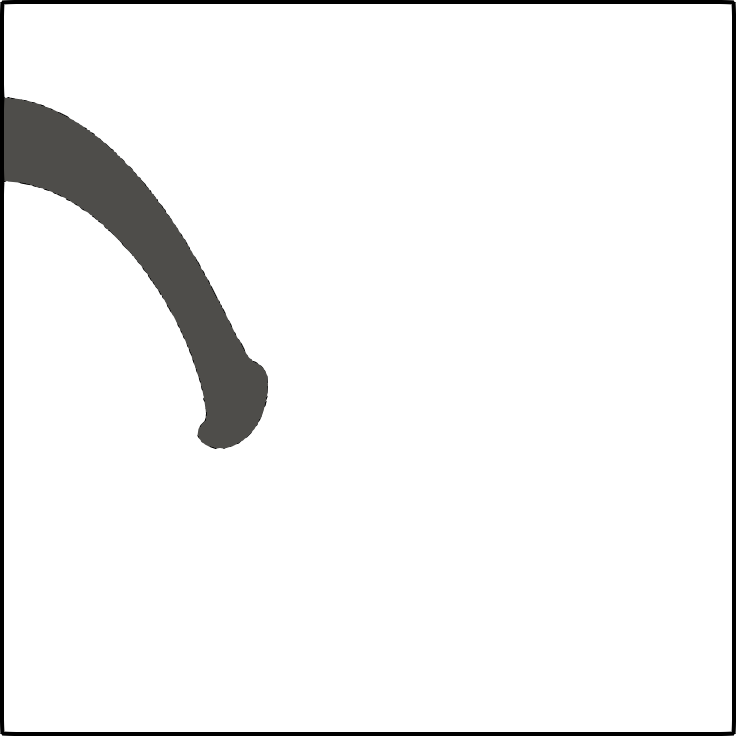} &
      \includegraphics[scale=0.13]{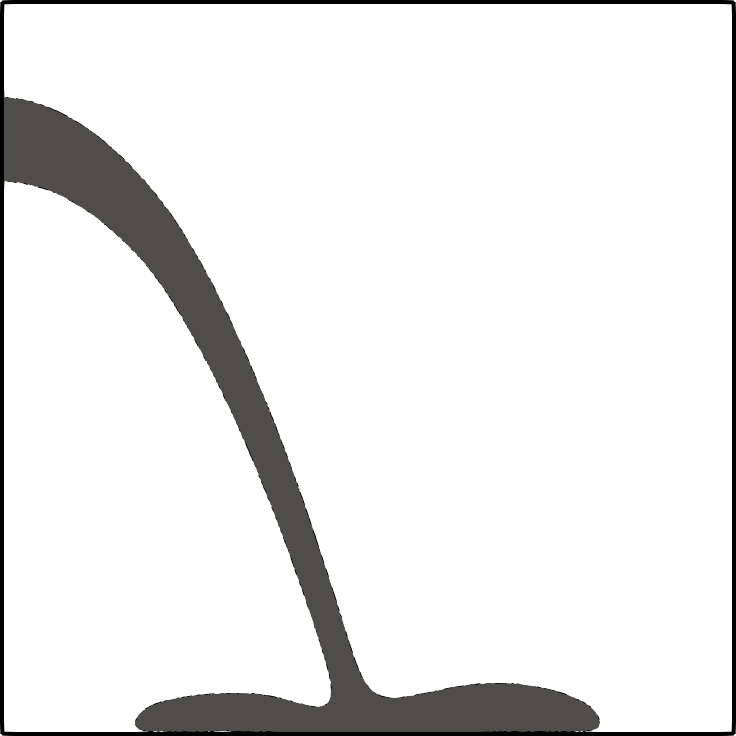} &
      \includegraphics[scale=0.13]{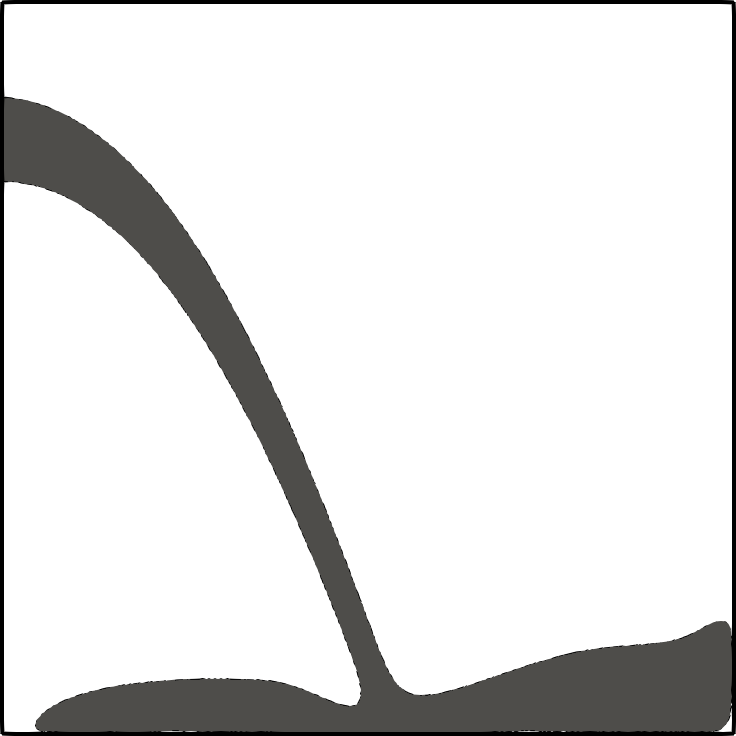} &
      \includegraphics[scale=0.13]{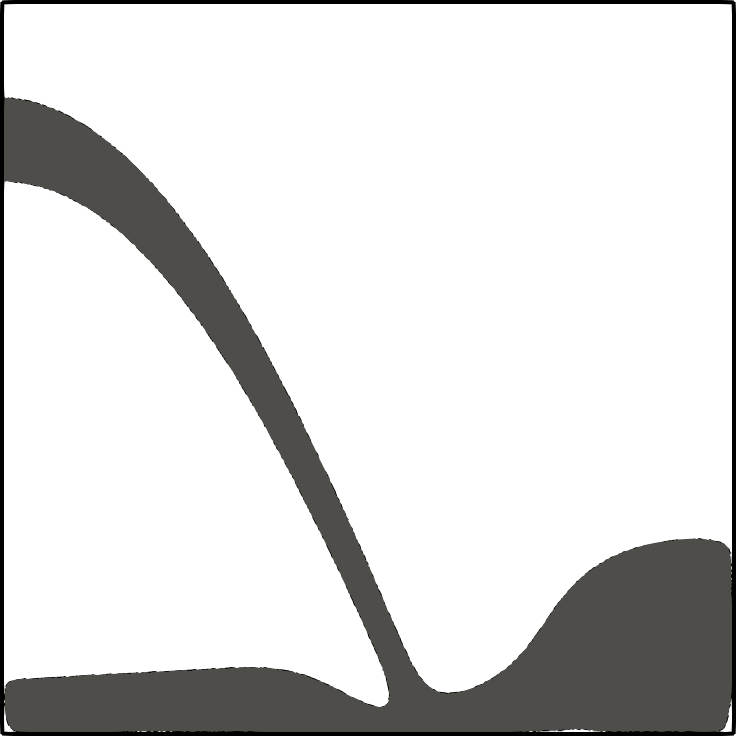} \\
      \includegraphics[scale=0.13]{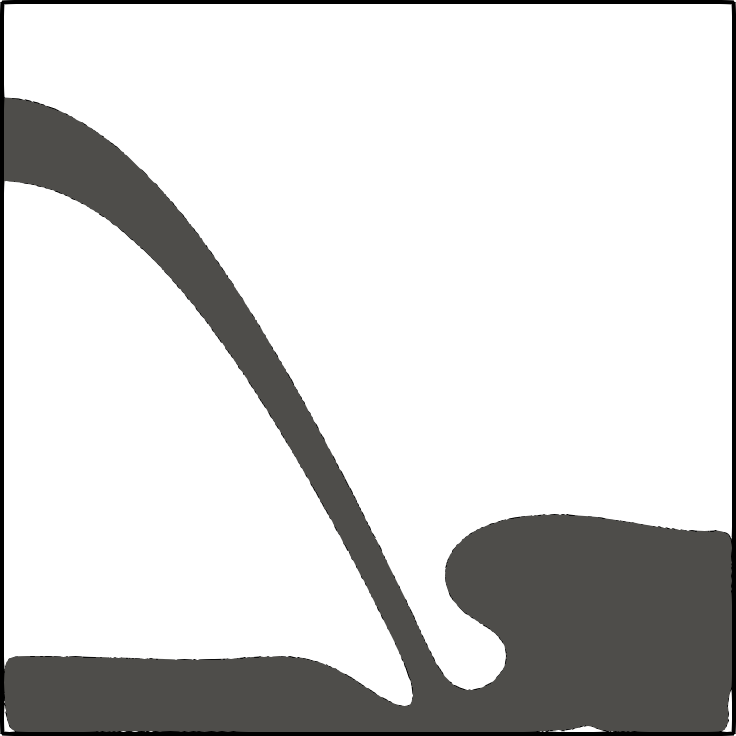} &
      \includegraphics[scale=0.13]{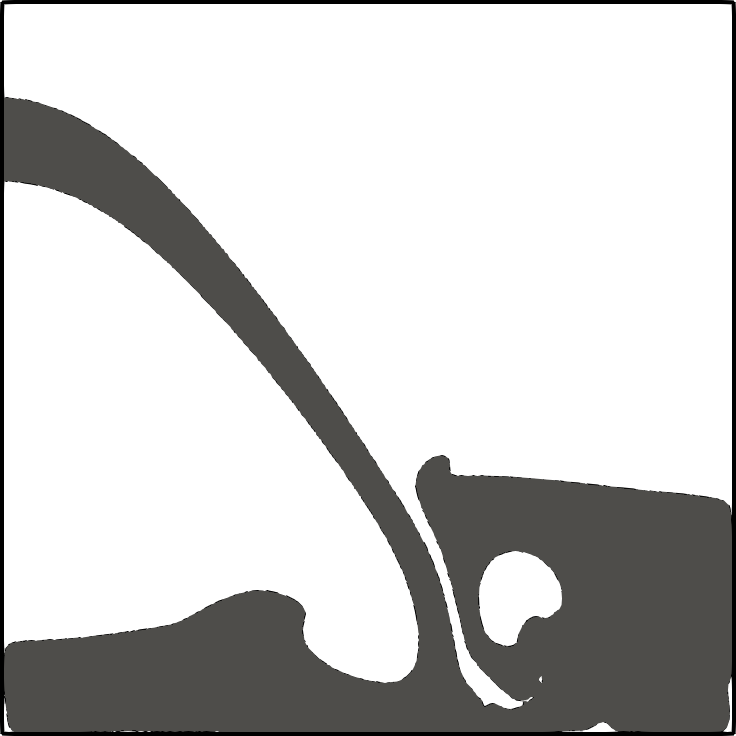} & 
      \includegraphics[scale=0.13]{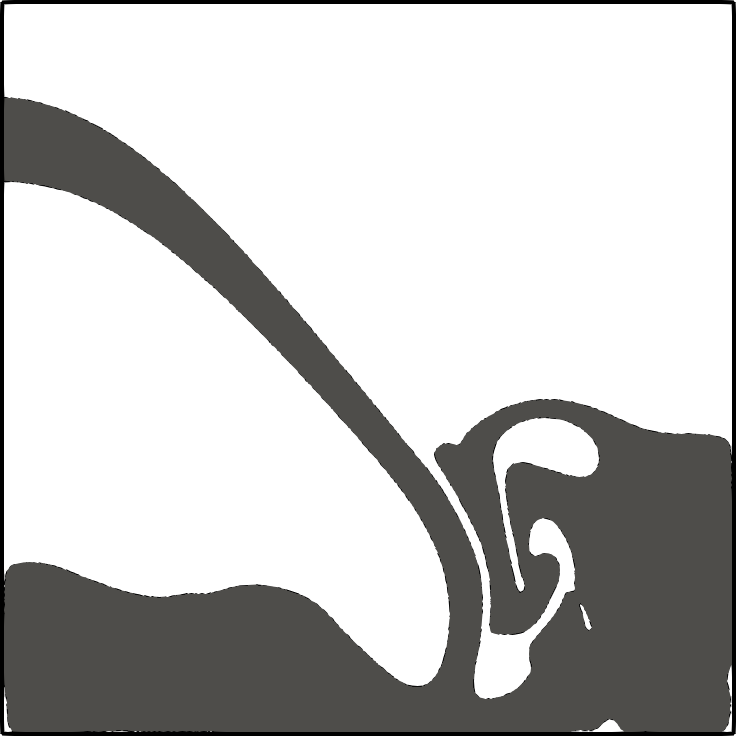} &
      \includegraphics[scale=0.13]{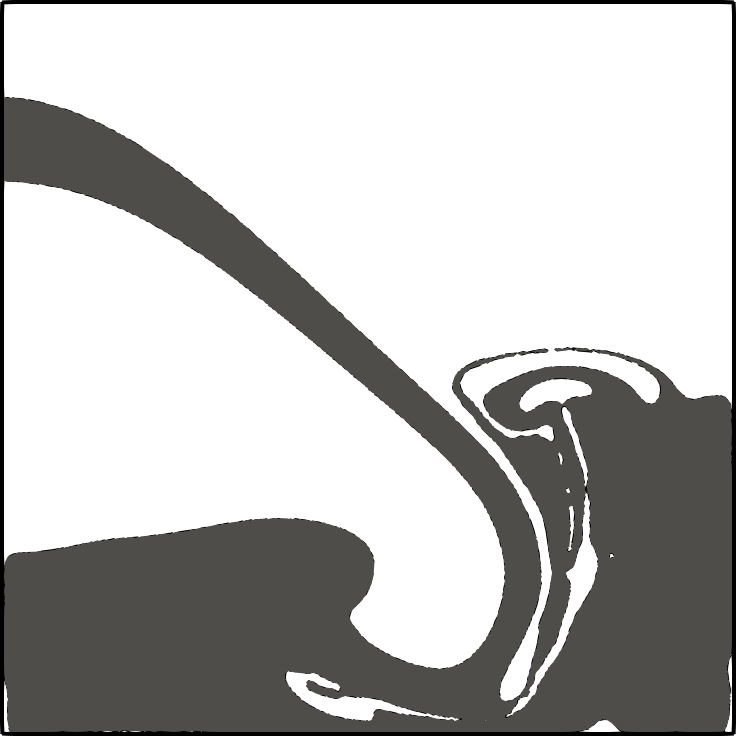}
  \end{tabular}}
  
  \subfloat[Test case 2]{
    \begin{tabular}{cccc} 
      \includegraphics[scale=0.13]{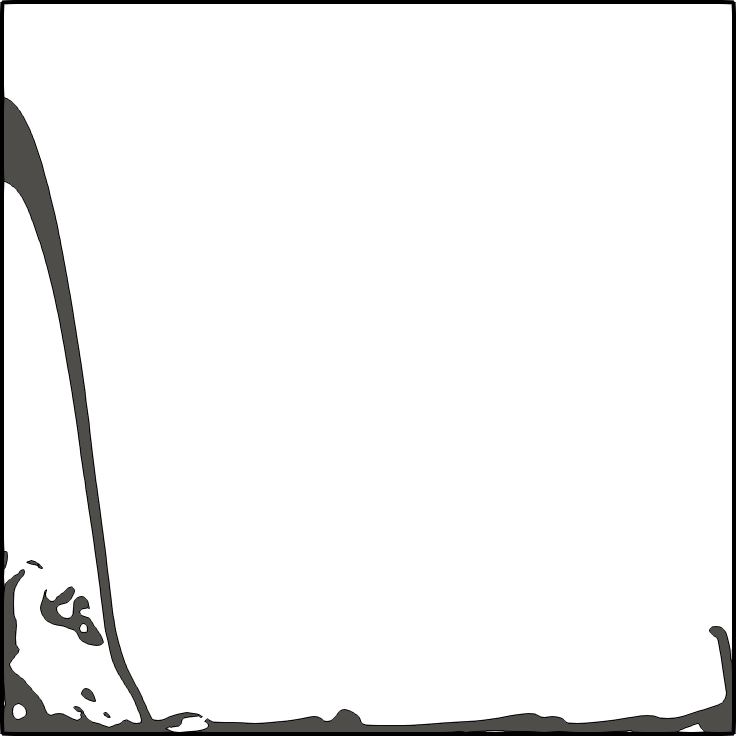} &
      \includegraphics[scale=0.13]{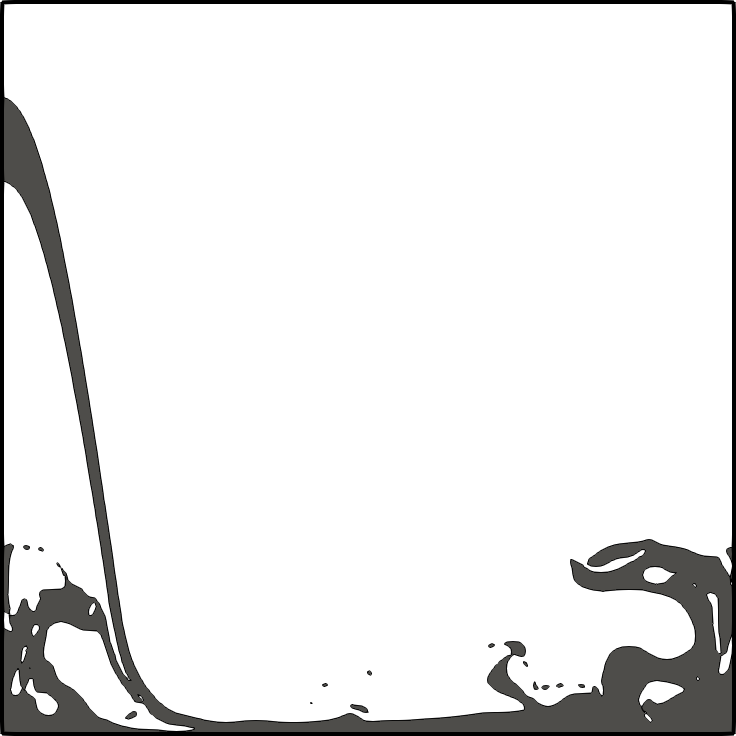} &
      \includegraphics[scale=0.13]{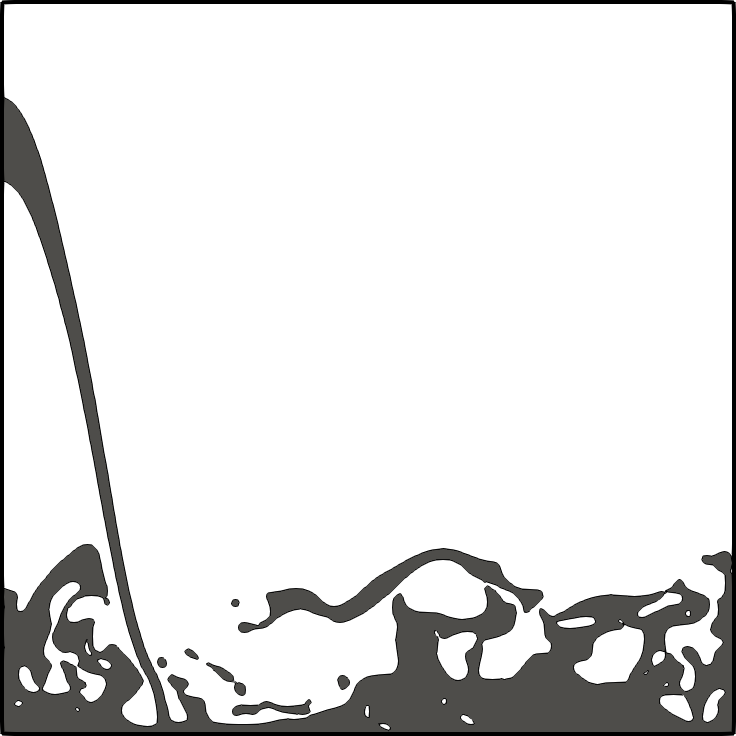} &
      \includegraphics[scale=0.13]{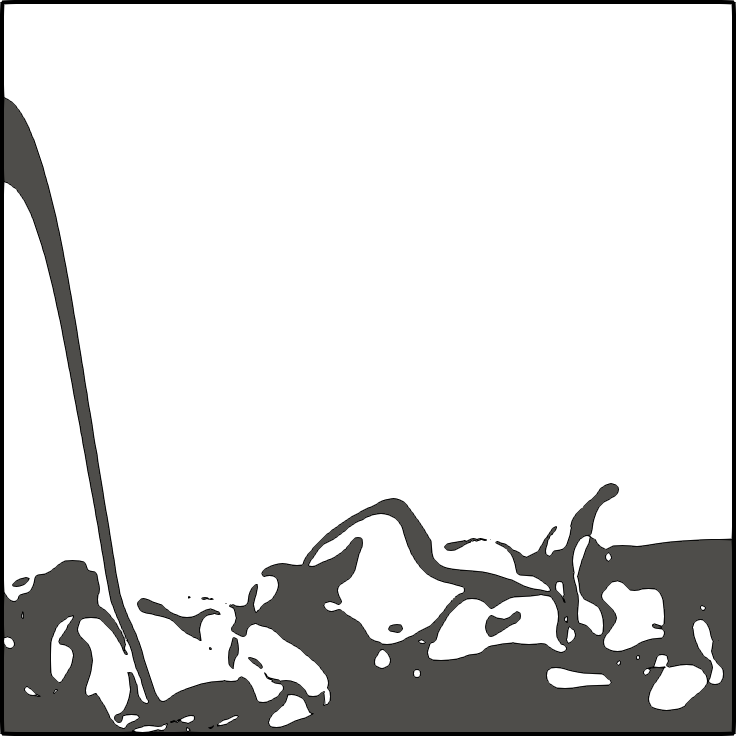} \\
      \includegraphics[scale=0.13]{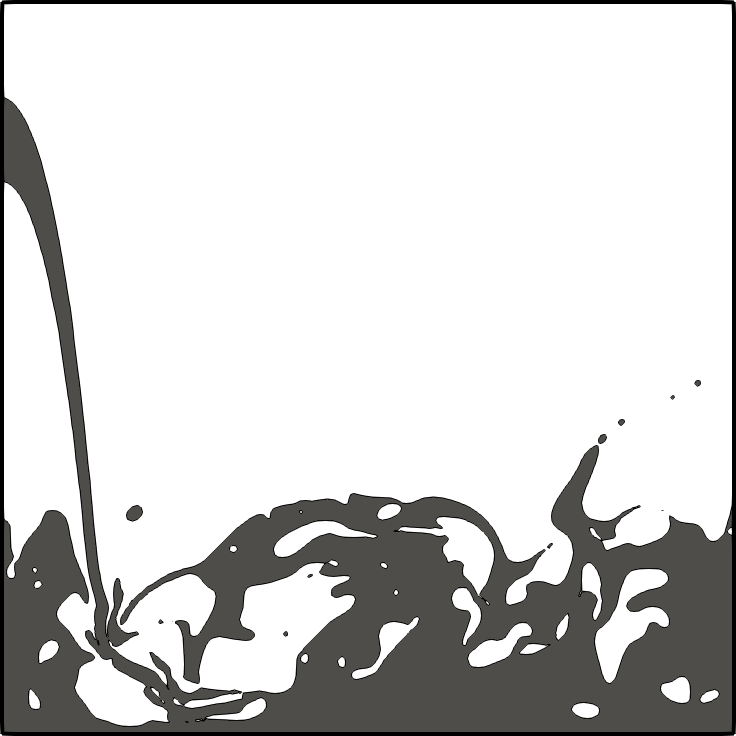} &
      \includegraphics[scale=0.13]{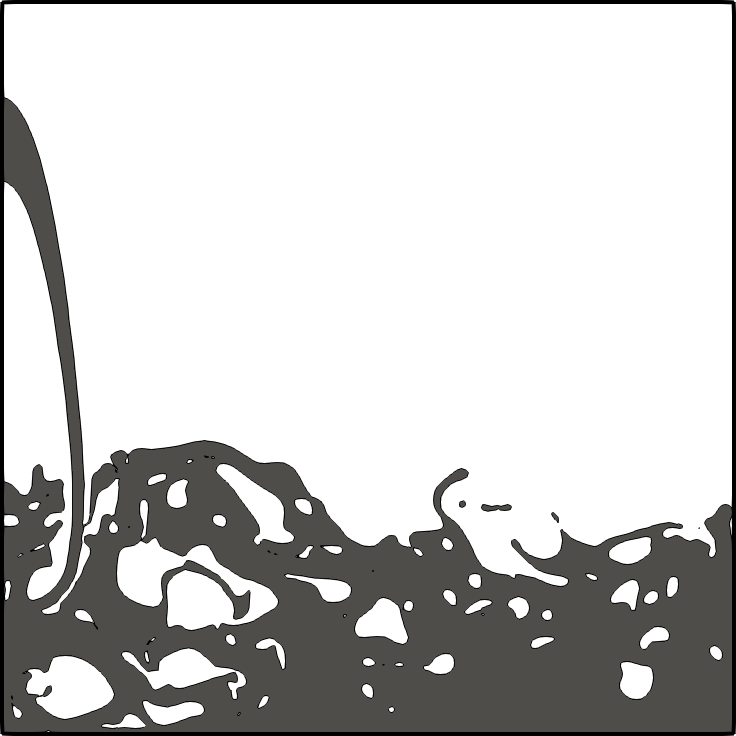} & 
      \includegraphics[scale=0.13]{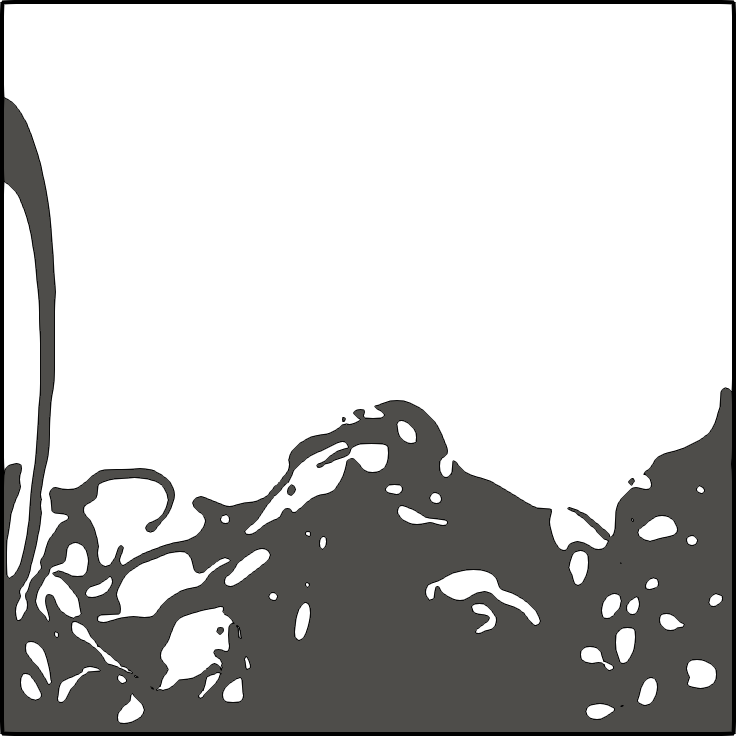} &
      \includegraphics[scale=0.13]{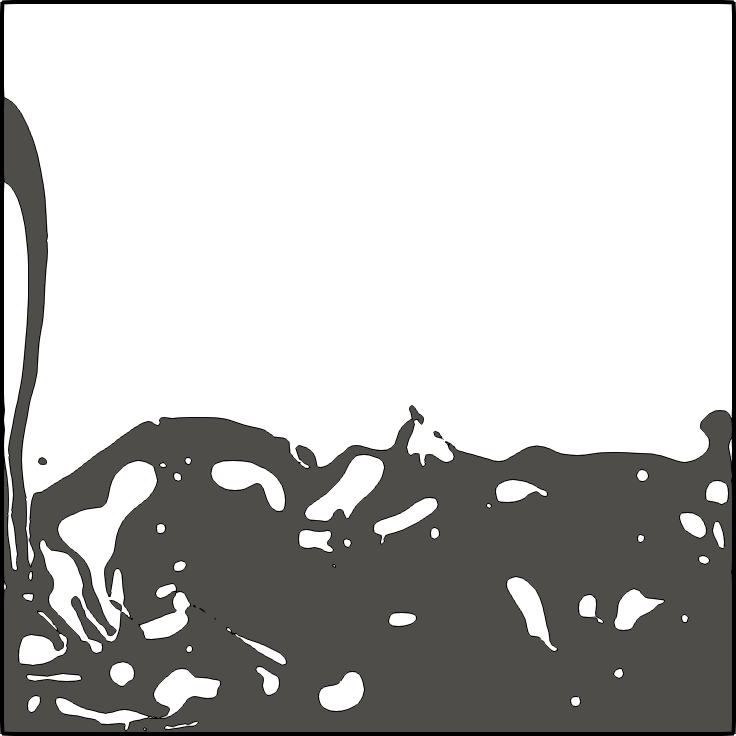}
    \end{tabular}}  
  \caption{Two-dimensional filling of a tank. For both test cases we show the heavy fluid phase
    at $t=0.3, ~0.6, ~0.9, ~1.2, ~1.5, ~1.8, ~2.1$ and $2.4$.}
  \label{fig:filling}
\end{figure}

In this section we consider the problem of filling an empty tank with a water-like fluid.
For this problem, the domain of interest is given by $\Omega=(0, 0.4)^2$.
The initial data consists of water in $W=\{(x,y)\in\Omega \st x<0.01, ~|y-0.325|<0.025\}$
and air in $A=\Omega\setminus W$. Both phases start at rest. 
The boundary $I=\{(x,y)\in\Omega \st x=0, ~|y-0.325|<0.025\}$ is set to inflow boundary.
At $I$ we set strongly $\bfu=(0.25,0)$. The rest of the left boundary is defined as slip boundary,
the bottom and the right boundaries are non-slip and the top is left open.
We consider two sets of material parameters. 
First we follow \cite[\S 11.3]{guermond2017conservative}
and reproduce qualitatively their results. In this case the material parameters are
\begin{align*}
  \rho_W=1000, \qquad
  \rho_A=1, \qquad
  \mu_W =1, \qquad
  \mu_A=1.8\times 10^{-2}, \qquad
  g=1, \qquad
  \sigma=0.
\end{align*}
Then we consider more realistic parameters given by \eqref{material_parameters},
without surface tension effects; i.e., $\sigma=0$.
We use an unstructured mesh with maximum element size of $h_e=2\times 10^{-3}$, which
corresponds to 126,888 triangular elements.
The water phase at different times are shown in figure \ref{fig:filling}.
Note the difference in the behavior of the solution due to different gravity and viscosity parameters. 

\subsection{Three-dimensional experiments}

\subsubsection{Dam break problem with obstacle}
\begin{figure}[!ht]
  \centering
  \includegraphics[scale=0.25]{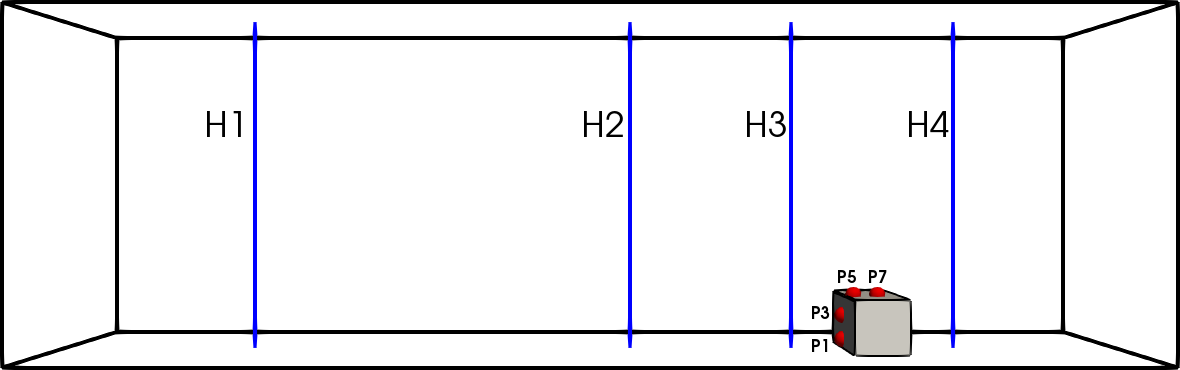}
  \caption{Water height and pressure gauges shown in blue and red, respectively.}
\label{fig:marin_location_gauges}
\end{figure}

\begin{figure}[!ht]
  \subfloat[Water height gauges. From left to right we show $H_1, ~H_2, ~H_3$ and $H_4$.]{
    \includegraphics[scale=0.18]{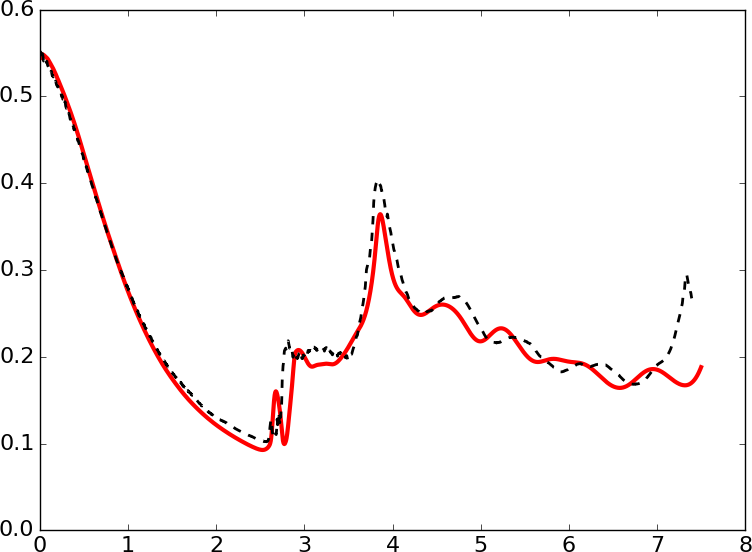}\qquad
    \includegraphics[scale=0.18]{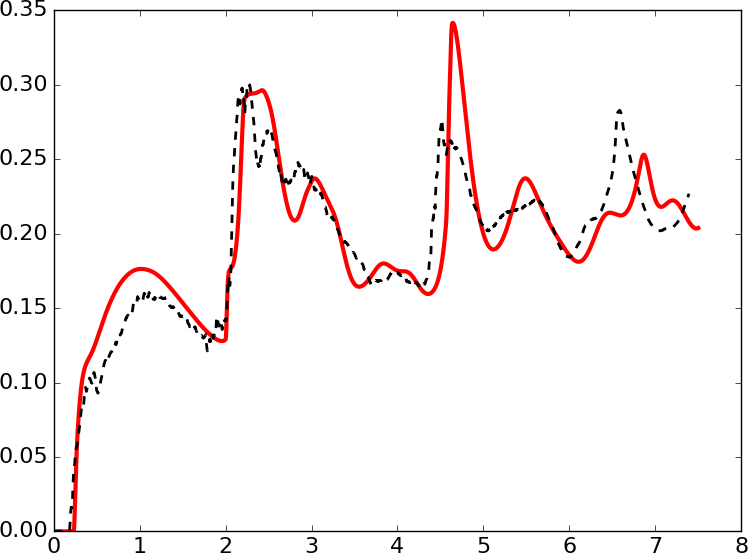}\qquad
    \includegraphics[scale=0.18]{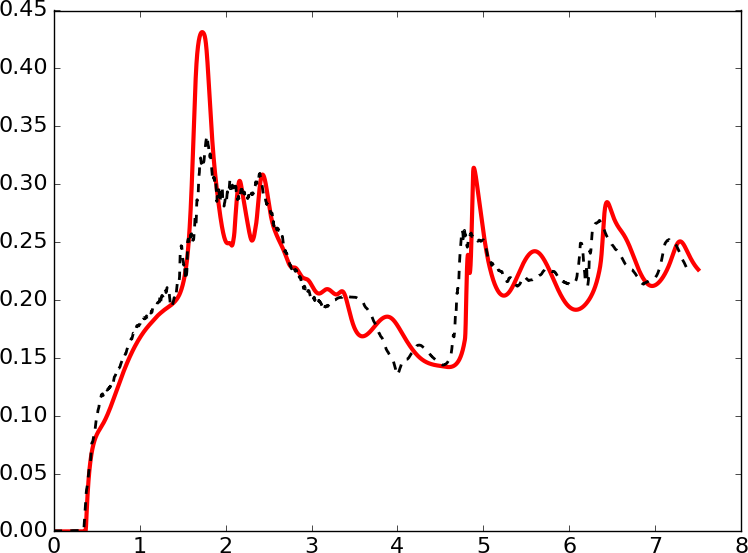}\qquad
    \includegraphics[scale=0.18]{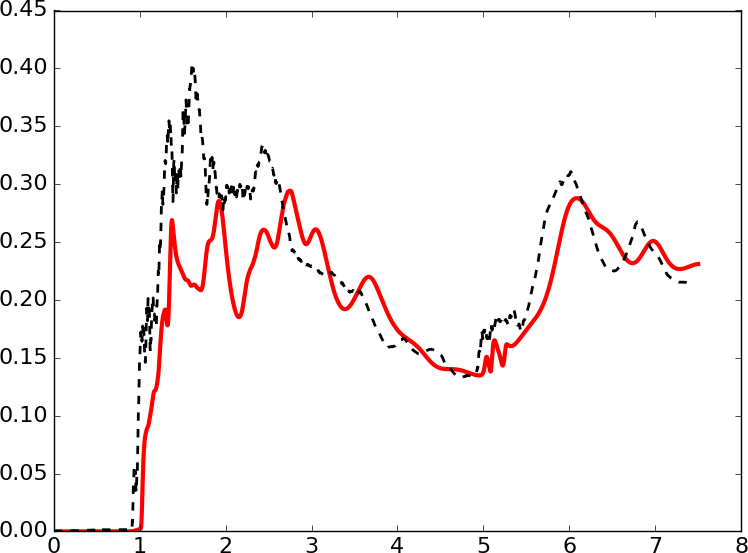}\qquad
  }
  
  \subfloat[Pressure gauges. From left to right we show $P_1, ~P_3, ~P_5$ and $P_7$.]{
    \includegraphics[scale=0.18]{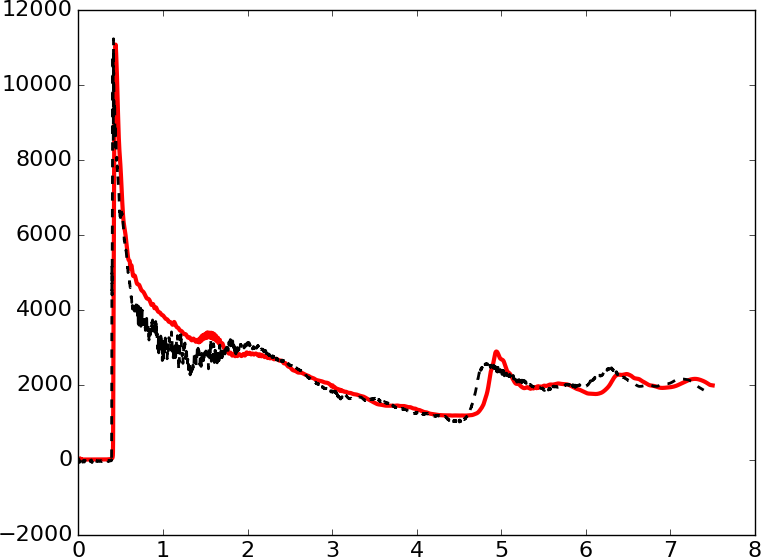}\qquad
    \includegraphics[scale=0.18]{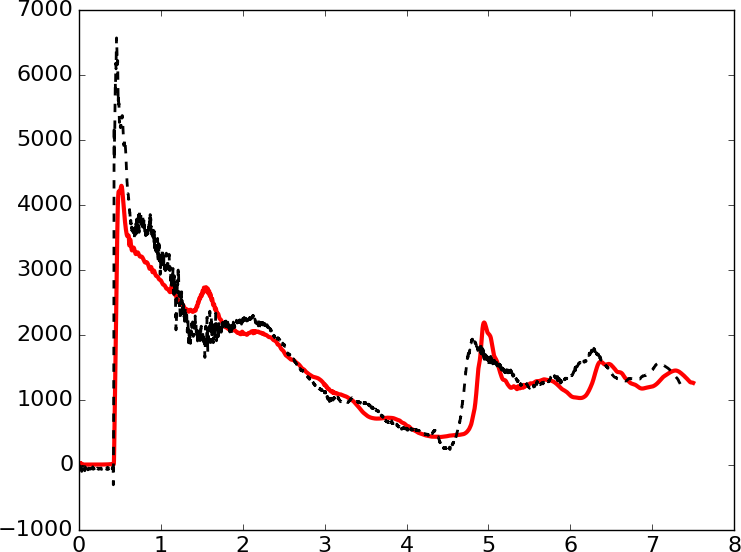}\qquad
    \includegraphics[scale=0.18]{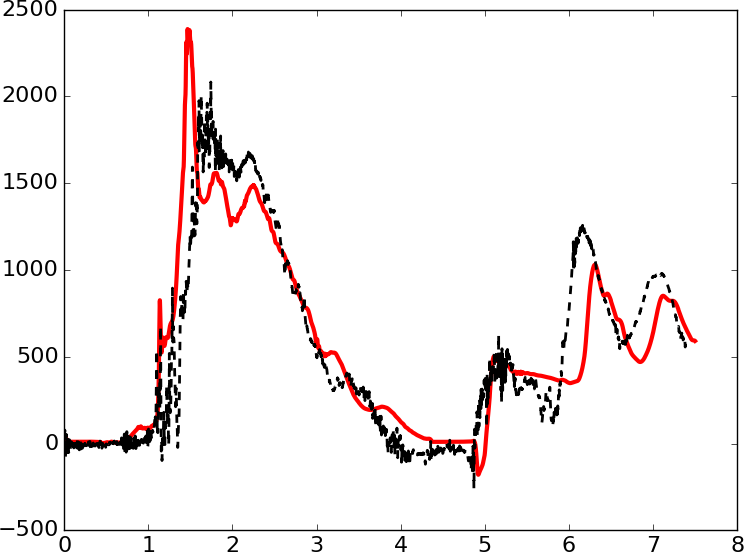}\qquad
    \includegraphics[scale=0.18]{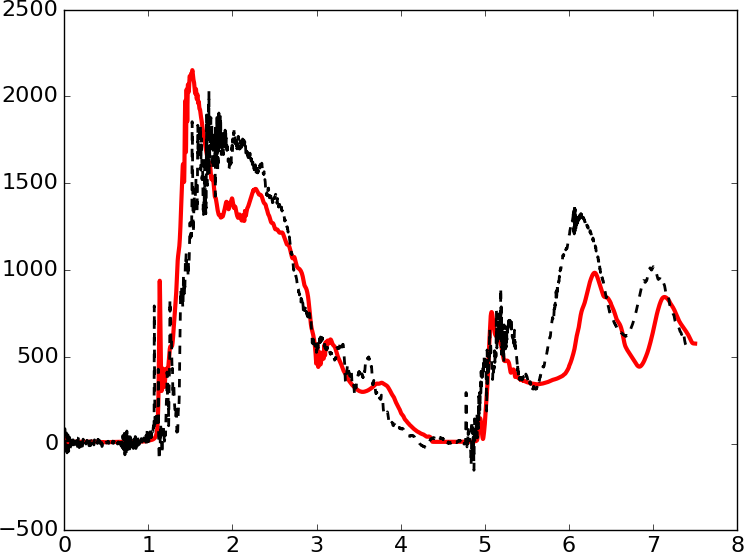}\qquad
  }    
  \caption{Water height and pressure gauges. The experimental data is plotted in dashed black
    and the numerical solution in solid red.}
\label{fig:marin_gauges}
\end{figure}

\begin{figure}[!ht]
  \centering
  \begin{tabular}{ccc} 
  \includegraphics[scale=0.16]{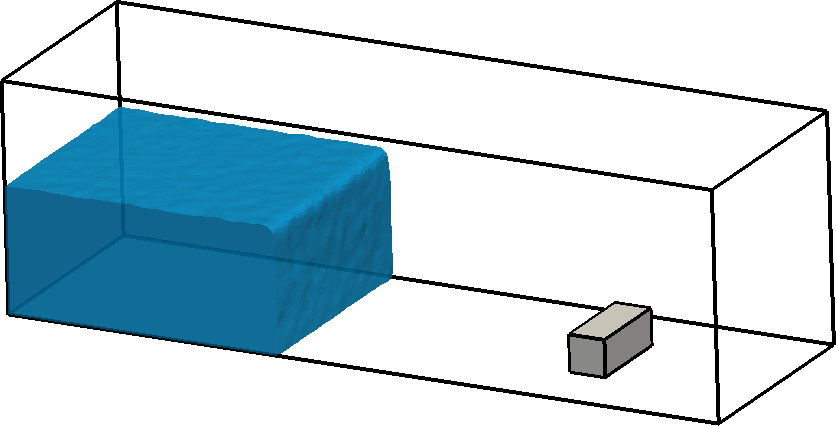} &
  \includegraphics[scale=0.16]{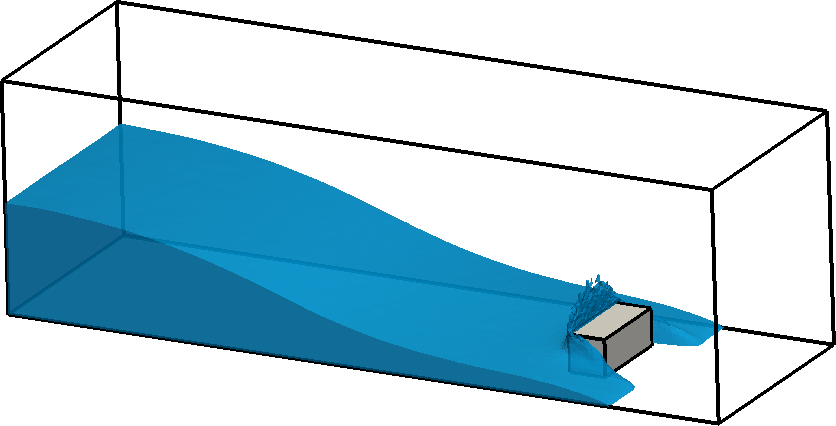} & 
  \includegraphics[scale=0.16]{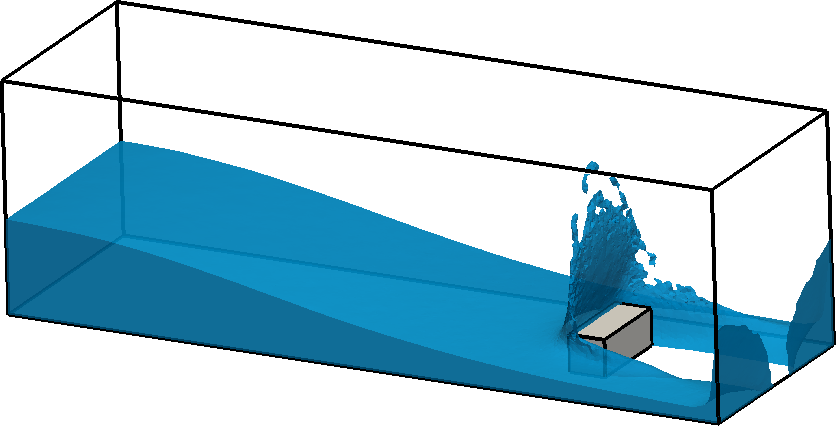} \\
  \includegraphics[scale=0.16]{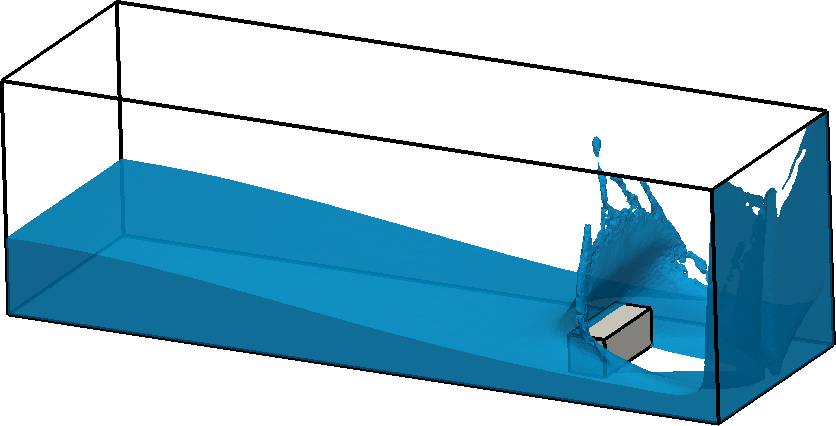} &
  \includegraphics[scale=0.16]{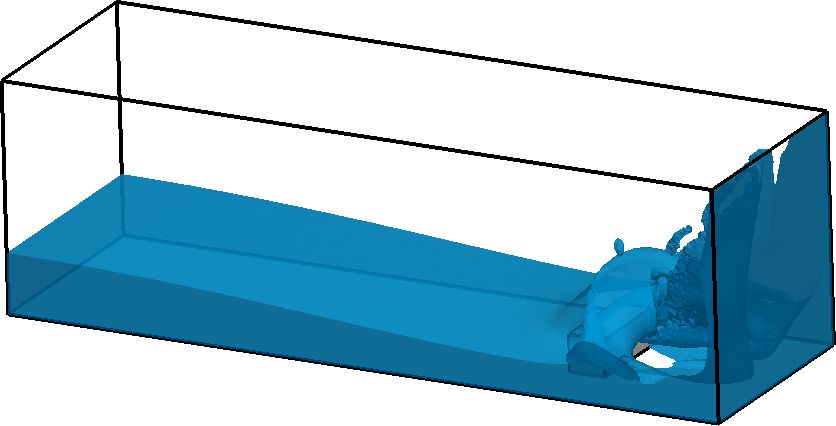} &
  \includegraphics[scale=0.16]{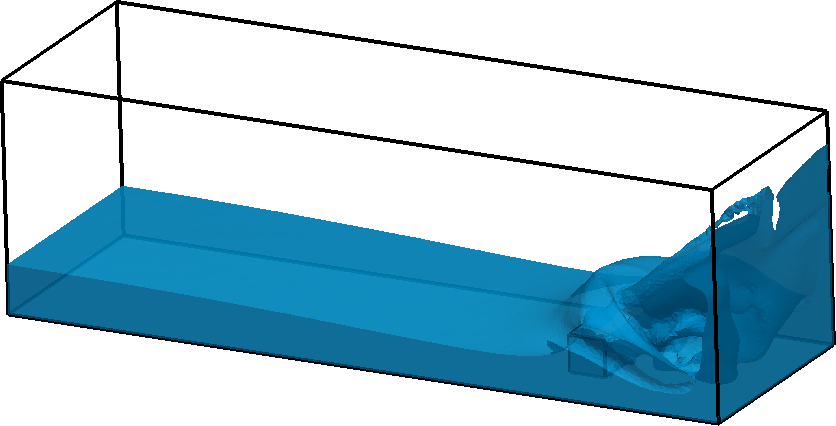} \\
  \includegraphics[scale=0.16]{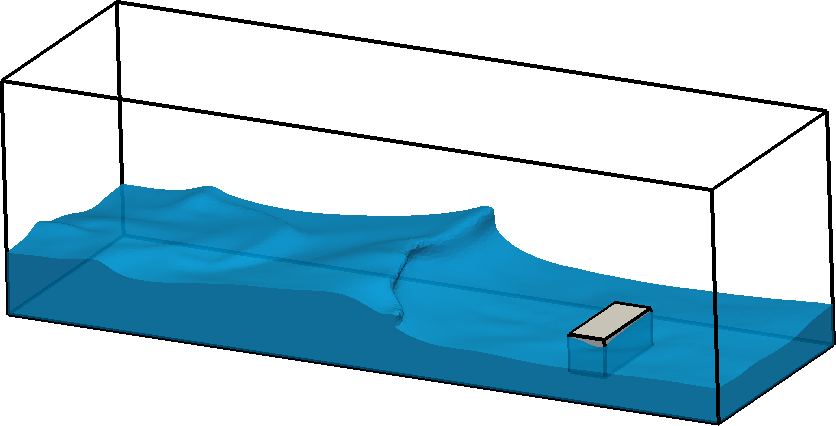} &
  \includegraphics[scale=0.16]{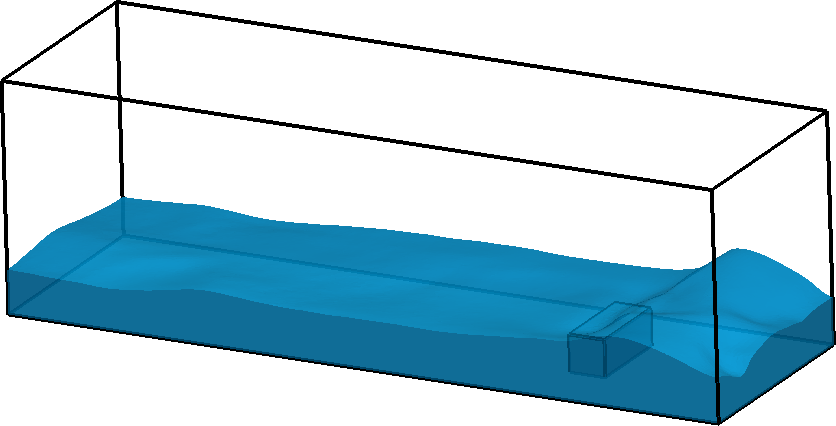} &
  \includegraphics[scale=0.16]{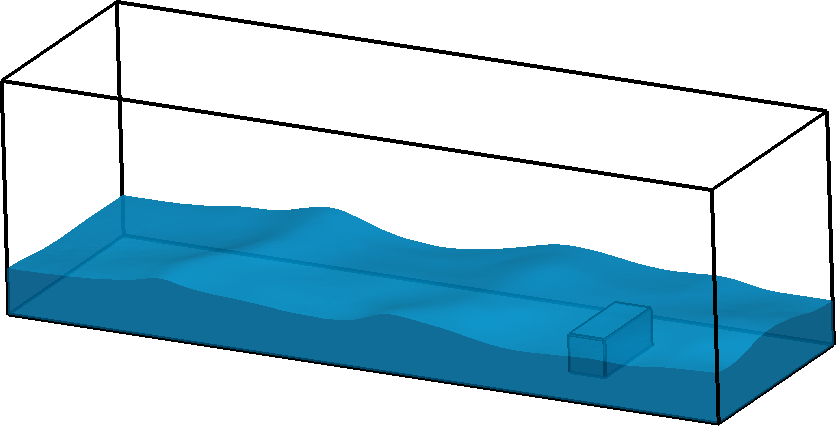}
  \end{tabular}  
  \caption{Three-dimensional dam break problem with obstacle.
    From top to bottom and left to right, we show the water phase
    at $t=0, ~0.5, ~0.7, ~0.9, ~1.1, ~1.4, ~4.5, ~6.0$ and $7.5$.}
\label{fig:marin}
\end{figure}

The first three dimensional problem that we consider is a dam break with an obstacle located downstream.
The computational domain is $\Omega=(0,3.22)\times(0,1)\times(0,1)\setminus O_1$,
where $O_1=[2.39,2.55]\times[0.3,0.7]\times[0,0.16]$ represents the obstacle. 
An initial column of water at rest is located in the domain
$W=\{\bfx\in\Omega \st x\leq 1.22, ~z\leq 0.55\}$.
The rest of the domain; i.e, $A=\Omega\setminus W$, is filled with air. 
At $t=0$ the column of water starts to fall because of the action of gravity.
Non-slip boundary conditions are imposed everywhere except for the top, which is left open. 
Experiments for this problem were performed by the Maritime Research Institute Netherlands
(MARIN), see \cite{kleefsman2005water, kleefsman2005volume}.
In particular, pressure and water height gauges were located in different
points of interest.
In this work we report the results for only four pressure gauges located at
$P_1=(2.389, 0.526, 0.025)$, $P_3=(2.389, 0.526, 0.099)$,
$P_5=(2.414, 0.474, 0.165)$ and $P_7=(2.487, 0.474, 0.165)$
and four water height gauges located at
$H_1=(0.582, 0.5, z), ~H_2=(1.732, 0.5, z), ~H_3=(2.228, 0.5, z)$ and
$H_4=(2.724, 0.5, z),  ~ \forall z\in[0,1]$.
See figure \ref{fig:marin_location_gauges}. 
We consider an unstructured mesh with characteristic size $h_e=0.025$,
which corresponds to 2,552,090 tetrahedral elements.
In figure \ref{fig:marin_gauges}, we compare the numerical results versus the
experimental measurements for the pressure and water height gauges and in
figure \ref{fig:marin} we show the water phase at different times.

\subsubsection{Moses flow}

\begin{figure}[ht]
  \centering
  \begin{tabular}{ccc} 
    \includegraphics[scale=0.19]{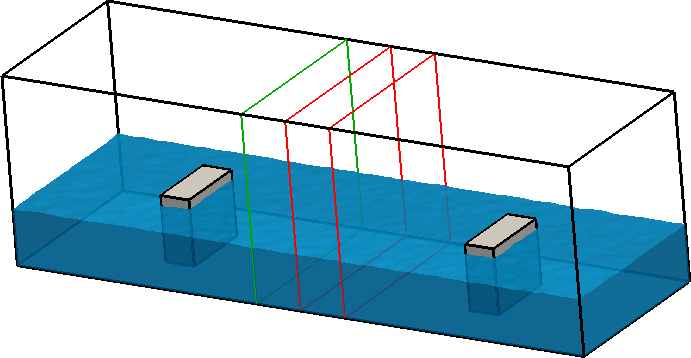} &
    \includegraphics[scale=0.19]{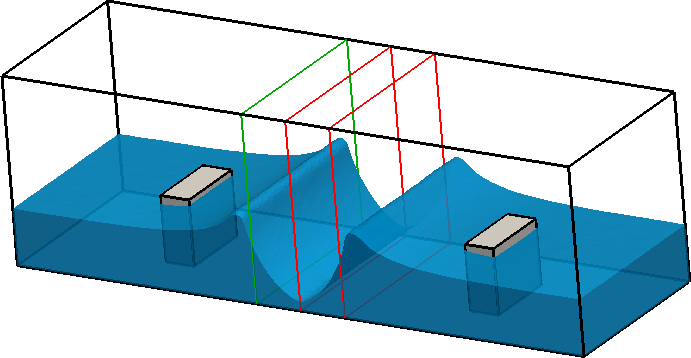} &
    \includegraphics[scale=0.19]{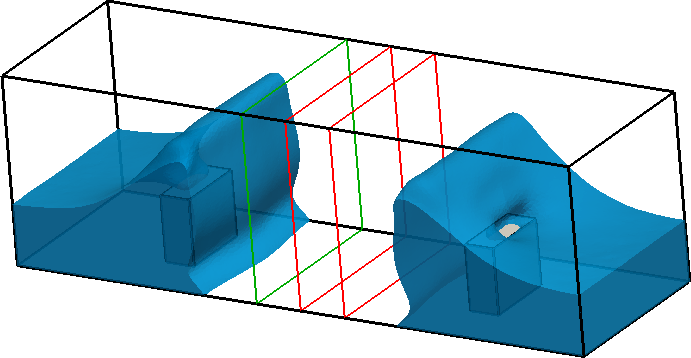} \\
    \includegraphics[scale=0.19]{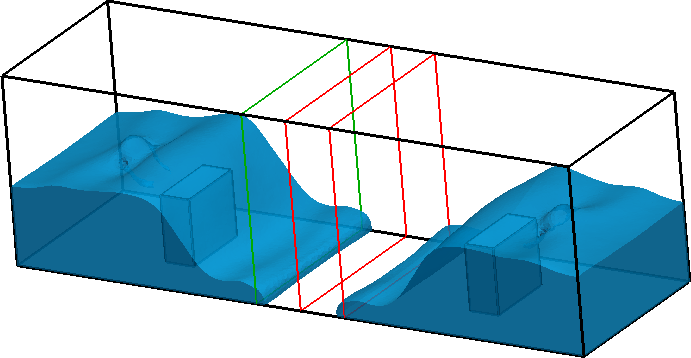} &
    \includegraphics[scale=0.19]{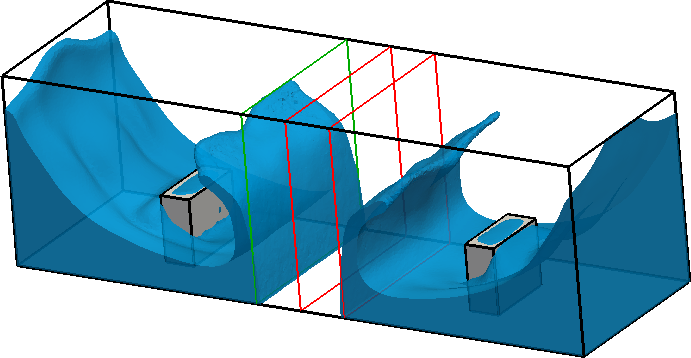} &
    \includegraphics[scale=0.19]{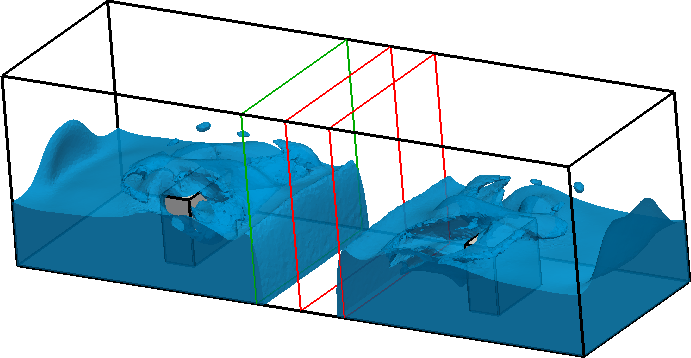} \\
    \includegraphics[scale=0.19]{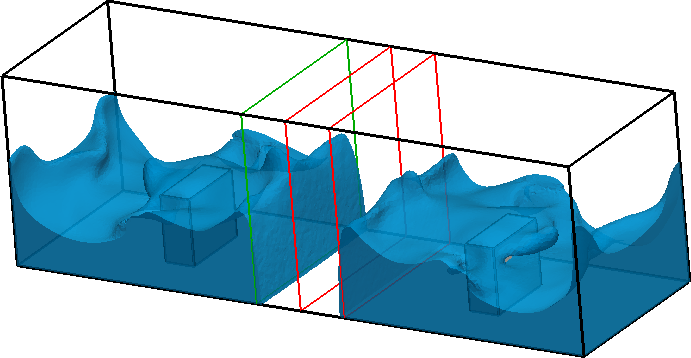} &
    \includegraphics[scale=0.19]{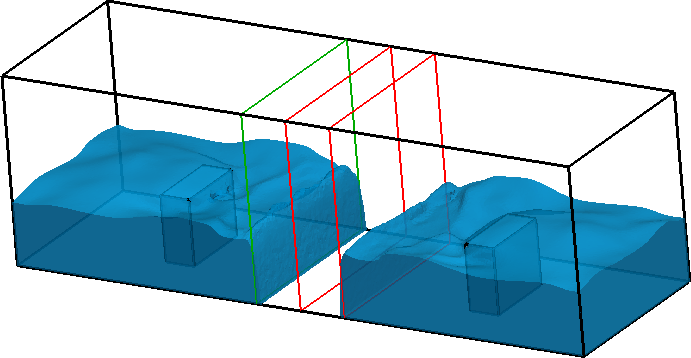} &
    \includegraphics[scale=0.19]{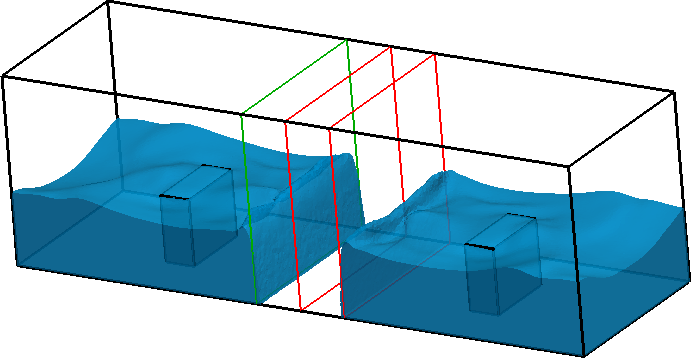}     
  \end{tabular}
  \caption{Three-dimensional Moses test.
    From top to bottom and left to right, we show the water phase
    at $t=0, ~0.1, ~0.4, ~0.6, ~0.9, ~1.3, ~2, ~3$, and $4$.}
  \label{fig:moses}
\end{figure}

In this simulation we consider the domain
$\Omega=(0,3.22)\times(0,1)\times(0,1)\setminus O_1\cup O_2$, where
$O_1=[0.6635,0.8245]\times[0.2985,0.7015]\times[0,0.35]$ and
$O_2=[2.3955,2.5565]\times[0.2985,0.7015]\times[0,0.35]$
represent two obstacles located at both ends of the domain.
Initially, the tank is partially filled with water and air at rest in the domains
$W=\{\bfx\in\Omega \st z\leq 0.3\}$ and $A=\Omega\setminus W$, respectively.
We apply an external force designed to split the water at the center of the tank.
Our aim with this problem is to test the robustness of the method with respect to
external (and aggressive) forces. 
The force is
\begin{align*}
  F(\bfx,t) =
  \begin{cases}
    -7.5\rho [1+\tanh(t)](10-9z), & \text{ if } 1.36 \leq x < 1.61, \\
    7.5\rho [1+\tanh(t)](10-9z), & \text{ if }  1.61 \leq x \leq 1.86
  \end{cases}.
\end{align*}
Note that the external force is scaled by the fluid density. 
All boundaries are considered to be non-slip except the top which is left open.
The grid is unstructured with element size $h_e=0.025$,
which corresponds to 2,522,647 tetrahedral elements.
In figure \ref{fig:moses} we show the water phase at different times. 
We remark that the solution appears to tend to a steady state. 

\section{Conclusions}\label{sec:conclusions}

In this work we presented a numerical model of incompressible, immiscible two-phase flow.
The general algorithm is driven by an operator splitting scheme in which a conservative level set method is first solved for a given
velocity field to evolve the interface. Afterwards, the 
interface (at the new location) is employed to obtain an updated velocity field.

We extended the conservative monolithic level set method presented in \cite{quezada2018monolithic} to make it more robust and to allow simpler initial conditions.
The original monolithic scheme contains a second-order term that penalizes deviations from the distance function, which was derived from parabolic and elliptic redistancing, see \cite{chan1999nonlinear, basting2014optimal}.
In our experience with two-phase flow simulations using the method in \cite{quezada2018monolithic}, we found it is not always possible to force the redistancing to emanate from the interface,
which results in incorrect solutions (of the viscous Eikonal equation \cite{mantegazza2003hamilton}). As this was not an issue with the earlier multistage conservative level set method presented in \cite{kees2011conservative} based on the classical Eikonal equation, we propose performing an extra pre-processing step to induce the redistancing to be propagated from the interface, effectively providing a better initial guess to the original monolithic scheme.
This extra step is done by solving the Eikonal equation directly.
It is important to remark that, thanks to the penalization embedded into the conservative level set,
solving the Eikonal equation in our case is not a computationally expensive task. This is true since the initial guess
in the non-linear iterative process already contains most of the features of the solution. 
We see this pre-processing step as a small but necessary improvement for the initial condition at each time step and consider a more robust monolithic scheme as future work.

We use a second order projection scheme proposed for for variable-density incompressible flows in \cite{guermond2009splitting} extended to solve the Navier-Stokes equations with variable density, viscosity, and time steps. For the applications
of consideration in this work, extra artificial viscosity is commonly needed to stabilize the momentum equations at high Reynolds number.
Special care is required when some of these techniques are used with unstructured meshes. 
In this work, we extend recent ideas on solving hyperbolic equations
via Discrete Maximum Principle (DMP) preserving continuous Galerkin finite elements. These methods are algebraic and suitable for both coarse and highly refined
unstructured meshes. We relax the DMP preserving features of these schemes, 
which allows simplification of the stabilization; however, in our numerical experiments
we always obtained well behaved solutions.
We remark that this stabilization technique is free of tunable parameters and capable of high-order accuracy with suitable changes to the underlying approximation spaces.

Once the individual methods are described, we present an extensive set of numerical examples in two- and three-dimensions.
In all these problems, we use the same numerical parameters which, in our opinion, demonstrates the robustness of the method. 
Finally, we provide an open source and freely downloadable computational framework to test and try the numerical examples 
presented in this work and others that the reader might be interested in. 

\section{Acknowledgments}

The work of Manuel Quezada de Luna was supported primarily by an
appointment to the Postgraduate Research Participation Program at the
U.S. Army Engineer Research and Development Center, Coastal and
Hydraulics Laboratory (ERDC-CHL) administrated by the Oak Ridge
Institute for Science and Education through an inter-agency agreement
between the U.S. Department of Energy and ERDC. Haydel Collins and
Chris Kees were supported by the ERDC University program and the ERDC
Future Investment Fund.
Permission was granted by the Chief of Engineers,
US Army Corps of Engineers, to publish this information.

\bibliographystyle{abbrvnat}
\bibliography{References.bib}

\end{document}